\documentclass[12pt, reqno]{amsart}

\usepackage[margin=.9in, a4paper, foot=.3in]{geometry}

\usepackage{graphicx, color}

\definecolor{darkblue}{rgb}{0,0,.5}
\definecolor{darkred}{rgb}{.5,0,0}
\definecolor{darkgreen}{rgb}{0,0.5,0}

\usepackage{hyperref}
\hypersetup {
    pdfstartview=FitH,
    colorlinks,
    citecolor=darkgreen,
    linkcolor=darkred,
    urlcolor=darkblue,
    linktoc=page
}

\let\oldtocsection=\tocsection
\let\oldtocsubsection=\tocsubsection
\let\oldtocsubsubsection=\tocsubsubsection

\renewcommand{\tocsection}[2]{\hspace{0em}\oldtocsection{#1}{#2}}
\renewcommand{\tocsubsection}[2]{\hspace{1em}\oldtocsubsection{#1}{#2}}
\renewcommand{\tocsubsubsection}[2]{\hspace{2em}\oldtocsubsubsection{#1}{#2}}

\usepackage{cite}

\usepackage[noBBpl,slantedGreek]{mathpazo}
\usepackage[mathcal]{euscript}

\usepackage{amssymb}
\usepackage{bm}
\allowdisplaybreaks
\numberwithin{equation}{section}

\usepackage{mathtools}

\usepackage{ifthen}

\newcommand {\svee}[2][\null]{\mathchoice     
  {#2%
  \ifthenelse{\equal{#1}{\null}}
    {\raisebox{-.15em}{\hbox{\large\char20}}}
    {\raisebox{-.15em}{\hbox{\large\char20}}\hskip -.3em_{#1}}}%
  {#2%
  \ifthenelse{\equal{#1}{\null}}
    {\raisebox{-.15em}{\hbox{\large\char20}}}
    {\raisebox{-.15em}{\hbox{\large\char20}}\hskip -.3em_{#1}}}%
  {#2%
  \ifthenelse{\equal{#1}{\null}}
    {\raisebox{-.15em}{\hbox{\char20}}}
    {\raisebox{-.15em}{\hbox{\char20}}\hskip -.3em_{#1}}}%     
  {#2%
  \ifthenelse{\equal{#1}{\null}}
    {\raisebox{-.15em}{\hbox{\small\char20}}}
    {\raisebox{-.15em}{\hbox{\small\char20}}\hskip -.3em_{#1}}}%   
}

\newcommand {\rmS}{\mathrm S}

\newcommand {\ch}{\mathrm{ch}}
\renewcommand {\det}{\mathrm{det}}
\newcommand {\diag}{\mathrm{diag}}
\renewcommand {\dim}{\mathrm{dim}}
\newcommand {\End}{\mathrm{End}}
\renewcommand {\exp}{\mathrm{exp}}
\newcommand {\id}{\mathrm{id}}

\newcommand {\Osc}{\mathrm{Osc}}
\newcommand {\sgn}{\mathrm{sgn}}
\newcommand {\tr}{\mathrm{tr}}

\newcommand {\bbC}{\mathbb C}
\newcommand {\bbE}{\mathbb E}
\newcommand {\bbZ}{\mathbb Z}

\newcommand {\calC}{\mathcal C}
\newcommand {\calL}{\mathcal L}
\newcommand {\calO}{\mathcal O}
\newcommand {\calQ}{\mathcal Q}
\newcommand {\calR}{\mathcal R}
\newcommand {\calT}{\mathcal T}
\newcommand {\calX}{\mathcal X}
\newcommand {\calY}{\mathcal Y}

\newcommand {\gothh}{\mathfrak h}
\newcommand {\gothg}{\mathfrak g}
\newcommand {\gothk}{\mathfrak k}

\newcommand {\gllpo}{\mathfrak{gl}_{\, l + 1}}
\newcommand {\hgothh}{\widehat{\mathfrak h}}
\newcommand {\hlgothg}{\widehat{\mathcal L}(\gothg)}
\newcommand {\lgothg}{\mathcal L(\mathfrak g)}
\newcommand {\sllpo}{\mathfrak{sl}_{\, l + 1}}
\newcommand {\tgothh}{\widetilde{\mathfrak h}}
\newcommand {\tlgothg}{\widetilde{\mathcal L}(\mathfrak g)}
\newcommand {\uqgllpo}{\mathrm U_q(\mathfrak{gl}_{\, l + 1})}
\newcommand {\uqlbp}{\mathrm U_q(\mathcal L(\mathfrak b_+))}
\newcommand {\uqlbm}{\mathrm U_q(\mathcal L(\mathfrak b_-))}
\newcommand {\uqlg}{\mathrm U_q(\calL(\mathfrak g))}
\newcommand {\uqlsllpo}{\mathrm U_q(\mathcal L(\mathfrak{sl}_{\, l + 1}))}
\newcommand {\uqlslii}{\mathrm U_q(\mathcal L(\mathfrak{sl}_2))}
\newcommand {\uqlsliii}{\mathrm U_q(\mathcal L(\mathfrak{sl}_3))}

\title{$\ell$-weights and factorization of transfer operators}

\author[A. V. Razumov]{A.~V.~Razumov}
\address{Institute for High Energy Physics, NRC ``Kurchatov Institute", 142281 Protvino, Mos\-cow region, Russia}
\email{Alexander.Razumov@ihep.ru}

\begin{document}

\addtolength {\jot}{3pt}

\begin{abstract}
We analyze the $\ell$--weights of the evaluation and $q$-oscillator representations of the quantum loop algebras $\uqlsllpo$ for $l = 1$ and $l = 2$ and prove the factorization relations for the transfer operators of the associated quantum integrable systems.
\end{abstract}

\maketitle

\tableofcontents

\section{Introduction}

This paper is devoted to the investigation of quantum integrable systems associated with the quantum loop algebras $\uqlsllpo$. The central object of the quantum group approach is the universal $R$-matrix being an element of the tensor product of two copies of the quantum loop algebra. The integrability objects are constructed by choosing representations for the factors of that tensor product. The consistent application of the method based on the quantum group theory was initiated by Bazhanov, Lukyanov and Zamolodchikov \cite{BazLukZam96, BazLukZam97, BazLukZam99}. They considered the quantum version of KdV theory. Later on this method proved to be efficient for studying other quantum integrable models. Within the framework of this approach, $R$-operators \cite{KhoTol92, LevSoiStu93, ZhaGou94, BraGouZhaDel94, BraGouZha95, BooGoeKluNirRaz10, BooGoeKluNirRaz11}, monodromy operators and $L$-operators were constructed \cite{BazTsu08,  BooGoeKluNirRaz10, BooGoeKluNirRaz11, BooGoeKluNirRaz13, Raz13, BooGoeKluNirRaz14a}. The corresponding sets of functional relations were discovered \cite{BazHibKho02, Koj08, BazTsu08, BooGoeKluNirRaz14a, BooGoeKluNirRaz14b, NirRaz14}.\footnote{For the terminology used we refer to the paper \cite{BooGoeKluNirRaz14a} and section \ref{s:urm} of the present paper.} Recently the quantum group approach was used to derive and investigate equations satisfied by the reduced density operators of the quantum chains related to an arbitrary loop algebra \cite{KluNirRaz20, Raz20}. 

Functional relations satisfied by commuting integrability objects are known to be a powerful tool for solving quantum integrable models. Usually, they are obtained by using the appropriate fusion rules for representations of the quantum loop algebra \cite{Res83a, Res83, KulRes86, BazRes90, KluPea92, KunNakSuz94, KunNakSuz11}. 

In fact, the most important functional relation is the factorized representation of the transfer operator. All other relations appears to be its consequences, see the papers \cite{NirRaz16a, BooGoeKluNirRaz14b} for $l = 1$ and $l = 2$. To prove the factorization relations a direct operator approach was used in those papers. For the higher ranks, the computational difficulties that arise seem to be almost insurmountable. In this paper, we propose to use a different approach based on the analysis of the $\ell$-weights of the representations. We demonstrate the effectiveness of the method for the cases $l = 1$ and $l = 2$.

In section 2 we define the quantum group $\uqgllpo$ and discuss its Verma modules. These modules are used later to define in section 3 the evaluation representations of the quantum loop algebras $\uqlsllpo$. We give two definitions of $\uqlsllpo$, the first in terms of the Drinfeld--Jimbo generators, and the second is the second Drinfeld's realization. The need for two definitions is that the first is convenient for defining evaluation representations, and the generators used in the second definition contain an infinite commutative subalgebra used to define $\ell$-weights. In section 4 we describe the category $\calO$ of $\uqlsllpo$-modules, introduce the concept of $\ell$-characters and define the Grothendieck ring of $\calO$. In the same section we define the $q$-oscillator representations of the Borel subalgebra $\uqlbp$ of $\uqlsllpo$. They are most convenient for our purposes. In the paper \cite{FreHer15} the prefundamental representations introduced in the paper \cite{HerJim12} are used.  We use the $q$-oscillator representations to construct $Q$-opera\-tors. In section 5 we define various types of integrability objects and discuss their properties. The factorization of the transfer operators for $l = 1$ and $l = 2$ is proved in section 6.

We fix the deformation parameter $\hbar$ in such a way that $q = \exp \hbar$ is not a root of unity, and define $q$-numbers by the equation
\begin{equation*}
[m]_q = \frac{q^m - q^{- m}}{q - q^{-1}}, \qquad m \in \bbZ.
\end{equation*}
In the present paper, an algebra, if it is not a Lie algebra, is a unital associative algebra. All algebras and vector spaces are assumed to be complex. By a tuple $(s_i)_{i \in I}$ we mean a mapping from a finite ordered set $I$ to some set of objects $S$. When a tuple has only one component or is used as a multi-index we omit the parentheses in the notation.

\section{\texorpdfstring{Quantum group $\uqgllpo$}{Quantum group Uq(gll+1)}} \label{s:2}

\subsection{Definition}

We start with a short reminder of some basics facts on  the Lie algebras $\gllpo$ and $\sllpo$. The standard basis of the standard Cartan subalgebra $\gothk$ of the general linear Lie algebra $\gllpo$ consists of the $(l + 1) {\times} (l + 1)$ matrices $K_a$, $a = 1, \ldots, l + 1$, defined by the equation\footnote{We denote by $\bbE_{a b}$ the usual matrix units.} 
\begin{equation*}
K_a = \bbE_{a a}.
\end{equation*}
Denote by $\epsilon_a$, $a = 1, \ldots, l + 1$, the elements of the dual basis. Below we often identify an element $\mu \in \gothk^*$ with the $(l + 1)$-tuple formed by its components $\mu_a = \langle \mu, \, K_a \rangle$, $a = 1, \ldots, l + 1$, with respect to this basis.

There are $l$ simple roots
\begin{equation*}
\alpha_i = \epsilon_i - \epsilon_{i + 1}, \qquad i = 1, \ldots, l.
\end{equation*}
The full system $\Delta^+$ of positive roots is formed by the roots
\begin{equation*}
\alpha_{i j} = \sum_{k = i}^{j - 1} \alpha_k = \epsilon_i - \epsilon_j, \qquad  1 \le i < j \le l + 1.
\end{equation*}
%The root
%\begin{equation*}
%\theta = \alpha_{1, \, l + 1} = \sum_{i = 1}^l \alpha_i
%\end{equation*}
%is the highest root of $\gllpo$.
The system of negative roots is $\Delta_- = - \Delta_+$, and the full root system is $\Delta = \Delta_+ \cup \Delta_-$.

The special Lie algebra $\sllpo$ is a subalgebra of $\gllpo$ formed by the traceless matrices. The standard basis of the standard Cartan subalgebra $\gothh$ of $\sllpo$ is formed by the matrices
\begin{equation*}
H_i = K_i - K_{i + 1}, \qquad i = 1, \ldots l.
\end{equation*}
As the positive and negative roots we take the restriction of $\alpha_{i j}$ and $- \alpha_{i j}$ to $\gothh$. For the simple roots we have
\begin{equation*}
\langle \alpha_i, \, H_j \rangle = a_{j i}, \label{aha}
\end{equation*}
where
\begin{equation*}
a_{i j} = {} - \delta_{i, \, j + 1} + 2 \delta_{i j} - \delta_{i + 1, \, j}. \label{aij}
\end{equation*}
The matrix $A = (a_{i j})_{i, \, j = 1}^l$ is the Cartan matrix of the Lie algebra $\sllpo$.

We define the quantum group $\uqgllpo$ as an algebra generated by the elements
\begin{equation*}
E_i, \quad F_i, \quad i = 1, \ldots, l, \qquad q^X, \quad X \in \gothk, \label{uqglg}
\end{equation*}
satisfying the following defining relations
\begin{gather*}
q^0 = 1, \qquad q^{X_1} q^{X_2} = q^{X_1 + X_2}, \\
q^X E_i \, q^{-X} = q^{\langle \alpha_i, \, X \rangle} E_i, \qquad q^X F_i \, q^{-X} = q^{-\langle \alpha_i, \, X \rangle} F_i, \\
[E_i, \, F_j] = \delta_{ij} \, \frac{q^{K_i - K_{i+1}} - q^{- K_i + K_{i + 1}}}{q - q^{-1}}.
\end{gather*}
In addition, there are the Serre relations. However, in the present paper we do not use their explicit form and, therefore, do not present it here.

\subsection{\texorpdfstring{Poincar\'e--Birkhoff--Witt basis}{Poincare-Birkhoff-Witt basis}} \label{s:pbwb}

The abelian group
\begin{equation*}
Q = \bigoplus_{i = 1}^l \bbZ \, \alpha_i
\end{equation*}
is called the root lattice of $\gllpo$. Assuming that
\begin{equation*}
E_i \in \uqgllpo_{\alpha_i}, \qquad F_i \in \uqgllpo_{- \alpha_i}, \qquad q^X \in \uqgllpo_0
\end{equation*}
for any $i = 1, \ldots, l$ and $X \in \gothk$, we endow the algebra $\uqgllpo$ with a $Q$-gradation. An element $x$ of $\uqgllpo$ is called a root vector corresponding to a root $\gamma$ of $\gllpo$ if $x \in \uqgllpo_\gamma$. In particular, $E_i$ and $F_i$ are root vectors corresponding to the roots $\alpha_i$ and $- \alpha_i$. Following Jimbo \cite{Jim86a}, we introduce the elements $E_{\gamma}$ and $F_{\gamma}$, corresponding to all roots $\gamma \in \Delta$, with the help of the relations
\begin{equation*}
E_{\alpha_{i, \, i + 1}} = E_i, \qquad E_{\alpha_{i j}} = E_{\alpha_{i, \, j - 1}} \, E_{\alpha_{j - 1, \, j}} - q \, E_{\alpha_{j - 1, \, j}} \, E_{\alpha_{i, \, j - 1}}, \quad j - i > 1,
\end{equation*} 
and
\begin{equation*} 
F_{\alpha_{i, \, i + 1}} = F_i, \qquad F_{\alpha_{i j}} = F_{\alpha_{j - 1, \, j}} \, F_{\alpha_{i, \, j - 1}} - q^{-1} F_{\alpha_{i, \, j - 1}} \, F_{\alpha_{j - 1, \, j}}, \quad j - i > 1.
\end{equation*}
It is clear that the elements $E_{\alpha_{i j}}$ and $F_{\alpha_{i j}}$ are root vectors corresponding to the roots $\alpha_{i j}$ and $- \alpha_{i j}$ respectively. They are linearly independent, and together with the elements $q^X$, $X \in \gothk$ are called the Cartan--Weyl generators of $\uqgllpo$. One can demonstrate that the ordered monomials constructed from the Cartan--Weyl generators form a Poincar\'e--Birkhoff--Witt basis of the quantum group $\uqgllpo$. In this paper we choose the following ordering. First endow $\Delta_+$ with the colexicographical order. It means that $\alpha_{i j} \preccurlyeq \alpha_{m n}$ if $j < n$, or if $j = n$ and $i < m$. This is a normal ordering of roots in the sense of \cite{LezSav74, AshSmiTol79}, see also \cite{Tol89}. In fact, this normal ordering is also a normal ordering of the roots of $\sllpo$. Now we say that a monomial is ordered if it has the form
\begin{equation*}
F_{\alpha_{i_1 j_1}} \ldots F_{\alpha_{i_r j_r}} \, q^X \, E_{\alpha_{m_1 n_1}} \ldots E_{\alpha_{m_s n_s}},
\end{equation*}
where $\alpha_{i_1 j_1} \preccurlyeq \cdots \preccurlyeq \alpha_{i_r j_r}$, $\alpha_{m_1 n_1} \preccurlyeq \cdots \preccurlyeq \alpha_{m_s n_s}$ and $X$ is an arbitrary element of $\gothk$. It is demonstrated in the paper \cite{NirRaz17b} that  such monomials really form a basis of $\uqgllpo$, see also the paper \cite{Yam89} for the case of an alternative ordering.

\subsection{Verma modules} \label{s:vm}

We use the standard terminology of the representation theory. In particular, we say that a $\uqgllpo$-module $V$ is a weight module if
\begin{equation*}
V = \bigoplus_{\mu \in \gothk^*}  V_\mu,
\end{equation*}
where
\begin{equation*}
V_{\bm \mu} = \{v \in V \mid q^X v = q^{\langle \mu, \, X \rangle} v \mbox{ for any } X \in \gothk \}.
\end{equation*}
The space $V_\mu$ is called the weight space of weight $\mu$, and a nonzero element of $V_\mu$ is called a weight vector of weight $\mu$. We say that $\mu \in \gothk^*$ is a weight of $V$ if $V_\mu \ne \{0\}$. A weight $\uqgllpo$-module $V$ is called a highest weight module of highest weight $\mu$ if there exists a weight vector $v \in V$ of weight $\mu$ such that
\begin{equation*}
E_i \, v = 0, \quad i = 1, \ldots, l, \quad \mbox{and} \quad V = \uqgllpo \, v.
\end{equation*}
The vector with the above properties is unique up to a scalar factor. We call it the highest weight vector of $V$.

Given $\mu \in \gothk^*$, denote by $\widetilde V^\mu$ the corresponding Verma $\uqgllpo$-module. This is a highest weight module of highest weight $\mu$. We denote by $\widetilde \pi^\mu$ the representation of $\uqgllpo$ corresponding to $\widetilde V^\mu$. Denote by $\bm m$ the $l (l + 1)/2$-tuple of non-negative integers $m_{i j}$, $1 \le i < j \le l + 1$, arranged in the colexicographical order of $(i, \, j)$. More explicitly,
\begin{equation*}
{\bm m} = (m_{12}, \, m_{13}, \, m_{23}, \, \ldots, \, m_{1j}, \, \ldots, \, m_{j-1, \, j}, \, \ldots, \,
m_{1, \, l + 1}, \, \ldots, \, m_{l, \, l + 1}).
\end{equation*}
The vectors
\begin{equation}
v_{\bm m} = F_{12}^{m_{12}} \, F_{13}^{m_{13}} \, F_{23}^{m_{23}} \, 
\ldots F_{1, \, j}^{m_{1, \, j}} \ldots F_{j - 1, \, j}^{m_{j - 1, \, j}} \ldots 
F_{1, \, l + 1}^{m_{1, \, l + 1}} \ldots F_{l, \, l + 1}^{m_{l, \, l + 1}} \, v_{\bm 0}, \label{vml}
\end{equation}
where for consistency we denote the highest weight vector by $v_{\bm 0}$, form a basis of $\widetilde V^\mu$. The explicit relations describing the action of the generators $q^X$, $E_i$ and $F_i$ of the quantum group $\uqgllpo$ on a general basis vector $v_{\bm m}$ are obtained in the paper \cite{NirRaz17b}.

Note that $\widetilde V^\mu$ is an infinite-dimensional $\uqgllpo$-module. However, if $\mu_i - \mu_{i + 1} \in \bbZ_{\ge 0}$ for all $i = 1, \ldots, l$, there is a unique maximal submodule of $\widetilde V^\mu$, such that the respective quotient module is simple and finite-dimensional. We denote this $\uqgllpo$-module and the corresponding representation by $V^\mu$ and $\pi^\mu$, respectively. Note that any finite-dimensional $\uqgllpo$-module can be constructed in this way.

The weights $\omega_a \in \gothk^*$, $a = 1, \ldots, l + 1$, defined as
\begin{equation*}
\omega_a = \sum_{b = 1}^a \epsilon_b = (\underbracket[.5pt]{1, \, \ldots, 1}_a, \, \underbracket[.5pt]{0, \, \ldots, \, 0}_{l + 1 - a}),
\end{equation*}
correspond to finite-dimensional representations called fundamental ones. The restriction of $\omega_i$, $i = 1, \ldots, l$, to $\gothh$ are the fundamental weights of $\sllpo$, so that
\begin{equation*}
\langle \omega_i, \, H_j \rangle = \delta_{i j}
\end{equation*}
for all $i, j = 1, \ldots, l$.

\section{\texorpdfstring{Quantum loop algebra $\uqlg$}{Quantum loop algebra Uq(L(g))}}

\subsection{\texorpdfstring{Definition in terms of Drinfeld--Jimbo generators}{Definition in terms of Drinfeld-Jimbo generators}} \label{ss:dtjmg}

Let $\gothg$ be a complex finite-dimen\-sional simple Lie algebra of rank $l$, $\gothh$ a Cartan subalgebra of $\gothg$, and $\Delta$ the root system of $\gothg$ relative to $\gothh$, see, for example, the books \cite{Ser01, Hum80}. Fix a system of simple roots $\alpha_i$, $i = 1, \ldots, l$. The corresponding coroots $\svee[i]{\alpha}$ form a basis of $\gothh$, so that
\begin{equation*}
\gothh = \bigoplus_{i = 1}^l \bbC \, \svee[i]{\alpha}.
\end{equation*}
The Cartan matrix $A = (a_{i j})_{i, \, j = 1}^l$ of $\gothg$ is defined by the equation
\begin{equation*}
a_{i j} = \langle \alpha_j, \, \svee[i]{\alpha} \rangle.
\end{equation*}
Note that any Cartan matrix is symmetrizable. It means that there exists a diagonal matrix $D = \diag(d_1, \ldots, d_l)$ such that the matrix $D A$ is symmetric and $d_i$, $i = 1, \ldots, l$, are positive integers. Such a matrix is defined up to a nonzero scalar factor. We fix the integers $d_i$ assuming that they are relatively prime.

Following Kac \cite{Kac90}, we denote by $\lgothg$ the loop algebra of $\gothh$, by $\tlgothg$ its standard central extension by the one-dimensional center $\bbC \, k$, and by $\hlgothg$ the Lie algebra obtained from $\tlgothg$ by adding a natural derivation $d$. By definition
\begin{equation*}
\hlgothg = \lgothg \oplus \bbC \, k \oplus \bbC \, d, 
\end{equation*}
and one uses as the Cartan subalgebra of $\hlgothg$ the space
\begin{equation*}
\hgothh = \gothh \oplus \bbC \, k \oplus \bbC \, d.
\end{equation*}
The Lie algebra $\hlgothg$ is isomorphic to the affine algebra associated with the extended Cartan matrix of $\gothg$. Here the basis coroots are $\svee[i]{\alpha}$, $i = 1, \ldots, l$, and
\begin{equation*}
\svee[0]{\alpha} = k - \sum_{i = 1}^l \svee[i]{a} \, \svee[i]{\alpha}.
\end{equation*}
The integers $\svee[i]{a}$, $i = 1, \ldots, l$, together with $ \svee[0]{a} = 1$ are the dual Kac labels of the Dynkin diagram associated with the extended Cartan matrix of $\gothg$. Thus, we have
\begin{equation*}
\hgothh = \bigoplus_{i = 0}^l \bbC \, \svee[i]{\alpha} \oplus \bbC \, d.
\end{equation*}
To introduce the corresponding simple roots, we identify the space $\gothh^*$ with the subspace of $\widehat \gothh^*$ defined as
\begin{equation*}
\{\lambda \in \widehat \gothh^* \mid \langle \lambda, \, k \rangle = 0, 
\ \langle \lambda, \, d \rangle = 0 \},
\end{equation*}
and denote by $\delta$ the element of $\hgothh^*$ defined by the equations
\begin{equation*}
\langle \delta, \, \svee[i]{\alpha} \rangle = 0, \quad i = 0, 1, \ldots, l, \qquad \langle \delta, \, d \rangle = 1.
\end{equation*}
Then the simple roots are $\alpha_i$, $i = 1, \ldots, l$, and
\begin{equation*}
\alpha_0 = \delta - \theta,
\end{equation*}
where
\begin{equation*}
\theta = \sum_{i = 1}^l a_i \alpha_i
\end{equation*}
is the highest root of $\Delta$. The integers $a_i$, $i = 1, \ldots, l$, together with $ a_0 = 1$ are the Kac labels of the Dynkin diagram associated with the extended Cartan matrix of $\gothg$.

One can demonstrate that the equation
\begin{equation*}
a_{i j} = \langle \alpha_j, \, \svee[i]{\alpha} \rangle, \qquad i, j = 0, 1, \ldots, l,
\end{equation*}
gives the entries of the extended Cartan matrix $A^{(1)}$ of the Lie algebra $\gothg$. Complementing the integers $d_i$, $i = 1, \ldots, l$, with a suitable integer $d_0$, one can demonstrate that the matrix $A^{(1)}$ is symmetrizable. Note that for $\gothg = \sllpo$, $d_i = 1$ for all $i = 0, 1, \ldots, l$.

The system of positive roots of the affine algebra $\hlgothg$ is
\begin{equation*}
\widehat \Delta_+ = \{\gamma + n \delta \mid  \gamma \in \Delta_+, \ n \in \bbZ_{\ge 0} \} \cup \{n \delta \mid n \in \bbZ_{>0} \} \cup \{(\delta - \gamma) + n \delta \mid  \gamma \in \Delta_+, \ n \in \bbZ_{\ge 0}\},
\end{equation*}
where $\Delta_+$ is the system of positive roots of the Lie algebra $\gothg$. The system of negative roots $\widehat \Delta_-$ of $\hlgothg$ is $\widehat \Delta_- = - \widehat \Delta_+$, and the full system of roots is
\begin{equation*}
\widehat \Delta = \widehat \Delta_+ \cup \widehat \Delta_- 
= \{ \gamma + n \delta \mid \gamma \in \Delta, \ n \in \bbZ \} \cup \{n \delta \mid n \in \bbZ \setminus \{0\} \}.
\end{equation*}

It is convenient for our purposes to denote
\begin{equation*}
\tgothh = \gothh \oplus \bbC \, k = \bigoplus_{i = 0}^l \bbC \, \svee[i]{\alpha}.
\end{equation*}
It is easy to show that for any $\lambda \in \gothh^*$ there is a unique element $\widetilde \lambda \in \widetilde \gothh^*$ such that
\begin{equation*}
\langle \widetilde \lambda, \, k \rangle = 0, \qquad \langle \widetilde \lambda, \, h \rangle = \langle \lambda, \, h \rangle, \quad h \in \gothh.
\end{equation*}

For each $i = 0, 1, \ldots, l$ we set
\begin{equation*}
q_i = q^{d_i}. 
\end{equation*}
and assume that
\begin{equation*}
q^\nu = \exp (\hbar \nu)
\end{equation*}
for any $\nu \in \bbC$. We define the quantum loop algebra $\uqlg$ as an algebra generated by the elements $e_i$, $f_i$, $i = 0, 1, \ldots, l$, and $q^h$, $h \in \tgothh$, satisfying the relations
\begin{gather*}
q^{\nu \, k} = 1, \quad \nu \in \bbC, \qquad q^{h_1} q^{h_2} = q^{h_1 + h_2}, \\
q^h e_i \, q^{-h} = q^{\langle \alpha_i, \, h \rangle} e_i, \qquad q^h f_i \, q^{-h} = q^{- \langle \alpha_i, \, h \rangle} f_i, \\
[e_i, \, f_j] = \delta_{ij} \, \frac{q_i^{\svee[i]{\alpha}} - q_i^{-\svee[i]{\alpha}}}{q^{}_i - q_i^{-1}}.
\end{gather*}
for all $i = 0, 1, \ldots, l$. There are also the Serre relations, whose explicit form is not used in the present paper.

The quantum loop algebra $\uqlg$ is a Hopf algebra. The comultiplication $\Delta$, the antipode $S$, and the counit $\varepsilon$ are given by the relations
\begin{gather*}
\Delta(q^h) = q^h \otimes q^h, \qquad \Delta(e^{}_i) = e^{}_i \otimes 1 + q_i^{\svee[i]{\alpha}} \otimes e^{}_i, \qquad \Delta(f^{}_i) = f^{}_i \otimes q_i^{- \svee[i]{\alpha}} + 1 \otimes f^{}_i, \\
S(q^h) = q^{- h}, \qquad S(e^{}_i) = - q_i^{- \svee[i]{\alpha}} e^{}_i, \qquad S(f^{}_i) = - f^{}_i \, q_i^{\svee[i]{\alpha}}, \\
\varepsilon(q^h) = 1, \qquad \varepsilon(e^{}_i) = 0, \qquad \varepsilon(f^{}_i) = 0.
\end{gather*}
We do not use these relations in the present paper. They are given to fix our conventions.

To distinguish from the tensor product of mappings, we denote the tensor product of representations $\varphi$ and $\psi$ of $\uqlg$ as
\begin{equation*}
\varphi \otimes_\Delta \psi = (\varphi \otimes \psi) \circ \Delta
\end{equation*}
and, similarly, the tensor product of the corresponding $\uqlg$-modules $V$ and $W$ as $V \otimes_\Delta W$.

The abelian group
\begin{gather*}
\widehat Q = \bigoplus_{i = 0}^l \bbZ \alpha_i
\end{gather*}
is called the root lattice of $\hlgothg$. Assuming that
\begin{gather*}
e_i \in \uqlg_{\alpha_i}, \qquad f_i \in \uqlg_{- \alpha_i}, \qquad q^h \in \uqlg_0
\end{gather*}
for any $i = 0, 1, \ldots, l$ and $h \in \tgothh$, we endow the algebra $\uqlg$ with a $\widehat Q$-gradation. An element $x$ of $\uqlg$ is called a root vector corresponding to a root $\gamma \in \widehat \Delta$ if $x \in \uqlg_\gamma$. One can construct linearly independent root vectors corresponding to all roots from~$\widehat \Delta$, see, for example, the papers \cite{TolKho92, KhoTol92, KhoTol93, KhoTol94} or the papers \cite{Bec94a, Dam98} for an alternative approach.

\subsection{Drinfeld's second realization} \label{s:2dr}

The quantum loop algebra $\uqlg$ can be realized in a different way \cite{Dri87, Dri88} as a $\bbC$-algebra with generators $\xi^\pm_{i, \, m}$ with $i = 1, \ldots, l$ and $m \in \bbZ$, $q^h$ with $h \in \gothh$, and $\chi_{i, \, m}$ with $i = 1, \ldots, l$ and $m \in \bbZ \setminus \{0\}$. They satisfy the defining relations
%\begin{gather*}
%q^0 = 1, \qquad q^{x_1} q^{x_2} = q^{x_1 + x_2}, \\
%[q^x, \, \chi_{j, \, m}] = 0, \qquad [\chi^{\mathstrut}_{i, \, m}, \, \chi^{\mathstrut}_{j, \, n}] = 0, \\
%q^x \xi^\pm_{i, \, m} q^{- x} = q^{\pm \langle \alpha_i, \, x \rangle} \xi^\pm_{i, \, m}, \qquad [\chi^{\mathstrut}_{i, \, m}, \, \xi^\pm_{j, n}] = \pm \frac{1}{m} [m \, a_{i j}]^{\mathstrut}_q \, \xi^\pm_{j, \, m + n}, \\
%\xi^\pm_{i, \, m + 1} \xi^\pm_{j, \, n} - q^{\pm a_{i j}} \, \xi^\pm_{j, \, n} \, \xi^\pm_{i, \, m + 1} = q^{\pm a_{i j}} \, \xi^\pm_{i, \, m} \, \xi^\pm_{j, \, n + 1} - \xi^\pm_{j, \, n + 1} \xi^\pm_{i, \, m}, \\
%[\xi^+_{i, \, n}, \, \xi^-_{j, \, m}] = \delta_{i j} \, \frac{q^{h_i} \phi^+_{i, \, n + m} - q^{- h_i} \phi^-_{i, \, n + m}}{q - q^{-1}}
%\end{gather*}
\begin{gather*}
q^0 = 1, \qquad q^{h_1} q^{h_2} = q^{h_1 + h_2}, \\
[q^h, \, \chi_{j, \, m}] = 0, \qquad [\chi^{\mathstrut}_{i, \, m}, \, \chi^{\mathstrut}_{j, \, n}] = 0, \\
q^h \xi^\pm_{i, \, m} q^{- h} = q^{\pm \langle \alpha_i, \, h \rangle} \xi^\pm_{i, \, m}, \qquad [\chi^{\mathstrut}_{i, \, m}, \, \xi^\pm_{j, n}] = \pm \frac{1}{m} [m \, a_{i j}]_{q^i} \, \xi^\pm_{j, \, m + n}, \\
\xi^\pm_{i, \, m + 1} \xi^\pm_{j, \, n} - q^{\pm a_{i j}}_i \, \xi^\pm_{j, \, n} \, \xi^\pm_{i, \, m + 1} = q^{\pm a_{i j}}_i \, \xi^\pm_{i, \, m} \, \xi^\pm_{j, \, n + 1} - \xi^\pm_{j, \, n + 1} \xi^\pm_{i, \, m}, \\
[\xi^+_{i, \, m}, \, \xi^-_{j, \, -m}] = \delta_{i j} \, \frac{q^{\svee[i]{\alpha}}_i - q^{- \svee[i]{\alpha}}_i}{q^{}_i - q^{-1}_i}, \\
[\xi^+_{i, \, m}, \, \xi^-_{j, \, n}] = \delta_{i j} \, \frac{q^{\svee[i]{\alpha}}_i \phi^+_{i, \, m + n}}{q^{}_i - q^{-1}_i},\quad m + n > 0, \quad [\xi^+_{i, \, m}, \, \xi^-_{j, \, n}] = - \delta_{i j} \, \frac{q^{- \svee[i]{\alpha}}_i \phi^-_{i, \, n + m}}{q^{}_i - q^{-1}_i}, \quad m + n < 0,
\end{gather*}
and the Serre relations whose explicit form is not important for our consideration. The quantities $\phi^\pm_{i, \, \pm m}$, $i = 1, \ldots, l$, $m \in \bbZ_{>0}$, are given by the formal power series
\begin{equation*}
1 + \sum_{m = 1}^\infty \phi^\pm_{i, \, \pm m} u^{\pm m} = \exp \left( \pm (q - q^{-1}) \sum_{m = 1}^\infty \chi_{i, \, \pm m} u^{\pm m} \right).
\end{equation*}
Stress that we use the definition of $\phi^\pm_{i, \, n}$ slightly different from the commonly used.

\subsection{Spectral parameters}

In applications to the theory of quantum integrable systems, one usually considers families of representations of a quantum loop algebra and its subalgebras, parameterized by a complex parameter called a spectral parameter. We introduce a spectral parameter in the following way. Assume that a quantum loop algebra $\uqlg$ is $\bbZ$-graded,
\begin{equation*}
\uqlg = \bigoplus_{m \in \bbZ} \uqlg_m, \qquad \uqlg_m \, \uqlg_n \subset \uqlg_{m + n},
\end{equation*}
so that any element $x \in \uqlg$ can be uniquely represented as
\begin{equation*}
x = \sum_{m \in \bbZ} x_m, \qquad x_m \in \uqlg_m.
\end{equation*}
Given $\zeta \in \bbC^\times$, we define the grading automorphism $\Gamma_\zeta$ by the equation
\begin{equation*}
\Gamma_\zeta(x) = \sum_{m \in \bbZ} \zeta^m x_m.
\end{equation*}
Now, for any representation $\varphi$ of $\uqlg$ we define the family $\varphi_\zeta$ of representations as
\begin{equation*}
\varphi_\zeta = \varphi \circ \Gamma_\zeta.
\end{equation*}
If $V$ is the $\uqlg$-module corresponding to the representation $\varphi$, we denote by $V_\zeta$ the $\uqlg$-module corresponding to the representation $\varphi_\zeta$.

In the present paper we endow $\uqlg$ with a $\bbZ$-gradation assuming that
\begin{equation*}
q^h \in \uqlg_0, \qquad e_i \in \uqlg_{s_i}, \qquad f_i \in \uqlg_{-s_i},
\end{equation*}
where $s_i$ are arbitrary integers. We denote
\begin{equation*}
s = \sum_{i = 0}^l a_i s_i,
\end{equation*}
where $a_i$ are the Kac labels of the Dynkin diagram associated with the extended Cartan matrix of $\gothg$.

\section{\texorpdfstring{Highest $\ell$-weight representations}{Highest l-weight representations}}

\subsection{\texorpdfstring{$\ell$-weights of $\uqlg$-modules}{l-weights of UqL(g)-modules}}  \label{ss:lw}

A $\uqlg$-module $V$ is called a weight module if
\begin{equation}
V = \bigoplus_{\lambda \in \gothh^*}  V_\lambda, \label{vvl}
\end{equation}
where
\begin{equation*}
V_\lambda = \{v \in V | q^h v = q^{\langle \widetilde \lambda, \, h} \rangle v \mbox{ for any } h \in \tgothh \}.
\end{equation*}
The notation $\widetilde \lambda$ for $\lambda \in \gothh^*$ is explained in section \ref{ss:dtjmg}. The space $V_\lambda$ is called the weight space of weight $\lambda$, and a nonzero element of $V_\lambda$ is called a weight vector of weight $\lambda$. We say that $\lambda \in \gothh^*$ is a weight of $V$ if $V_\lambda \ne \{ 0 \}$.

A $\uqlg$-module $V$ is said to be in category $\calO$ if
\begin{itemize}
\item[(i)] $V$ is a weight module all of whose weight spaces are finite-dimensional;
\item[(ii)] there exists a finite number of elements $\mu_1, \ldots, \mu_s \in \gothh^*$ such that every weight of $V$ belongs to the set $\bigcup_{i = 1}^s D(\mu_i)$ , where $D(\mu) = \{\lambda \in \gothh^* \mid \lambda \leq \mu \}$ with $\leq$ being the usual partial order in $\gothh^*$, see, for example, the book \cite{Hum80}.
\end{itemize}

%\begin{additions}
%\begin{equation*}
%k_i = q^{d_i h_i} = q_i^{h_i}
%\end{equation*}
%If $v \in V_\lambda$
%\begin{equation*}
%k_i v = q_i^{\langle \widetilde \lambda, \, h_i \rangle} v = q_i^{\langle \lambda, \, h_i \rangle} v
%\end{equation*}
%\begin{equation*}
%\lambda = \sum_{j = 1}^l \lambda_j \omega_j
%\end{equation*}
%\begin{equation*}
%k_i v = q_i^{\lambda_i} v
%\end{equation*}
%$\lambda = \omega_j$
%\begin{equation*}
%k_i v = q^{\delta_{i j}}_i v
%\end{equation*}
%
%\end{additions}
%
%\begin{additions}
% 
%A $\uqlg$-module $V$ in the category $\calO$ is called a highest weight module of highest weight $\lambda \in \gothh^*$ if there exists a weight vector $v \in V$ of weight $\lambda$ such that
%\begin{equation*}
%e_i v = 0, \qquad i = 1, \ldots, l, \quad \mbox{and} \quad V = \uqlsllpo \, v.
%\end{equation*}
%The vector with the above properties is unique up to a scalar factor. It is called the highest weight vector of $V$.
%
%\end{additions}

Let a $\uqlg$-module $V$ be in category $\calO$. The algebra $\uqlg$ contains an infinite-dimensional commutative subalgebra generated by the elements $\phi^\pm_{i, \, \pm m}$, $i = 1, \ldots, l$, $m \in \bbZ_{>0}$ and $q^h$, $h \in \widetilde \gothh$. We can refine the weight decomposition (\ref{vvl}) in the following way. Let $\lambda$ be a weight of $V$. By definition, the space $V_\lambda$ is finite-dimensional. The restriction of the action of the elements $\phi^\pm_{i, \, \pm m}$ to $V_\lambda$ constitutes a countable set of pairwise commuting linear operators on $V_\lambda$. Hence, there is a basis of $V_\lambda$ which consists of eigenvectors and generalized eigenvectors of all those operators, see, for example, the book \cite{Lax07}. This leads to the following definitions.

An $\ell$-weight is a triple
\begin{equation*}
\bm \Lambda = (\lambda, \,  \bm \Lambda^+ \! , \, \bm \Lambda^-),
\end{equation*}
where $\lambda \in \gothh^*$, $\bm \Lambda^+ = (\Lambda^+_i(u))^l_{i = 1}$ and $\bm \Lambda^- = (\Lambda^-_i(u^{-1}))^l_{i = 1}$ are $l$-tuples of formal power series
\begin{equation*}
\Lambda^+_i(u) = 1 + \sum_{m \in \bbZ_{> 0}} \Lambda^+_{i, \, m} u^m \in \bbC[[u]], \qquad \Lambda^-_i(u^{-1}) = 1 + \sum_{m \in \bbZ_{> 0}} \Lambda^-_{i, \, - m} u^{- m} \in \bbC[[u^{-1}]].
\end{equation*}
We denote by $\gothh^*_\ell$ the set of $\ell$-weights.

Define a surjective homomorphism $\varpi \colon \gothh^*_\ell \to \gothh^*$ by the relation
\begin{equation*}
\varpi(\bm \Lambda) = \lambda
\end{equation*}
if $\bm \Lambda = (\lambda, \, \bm \Lambda^+ \! , \, \bm \Lambda^-)$. Now we have
\begin{equation*}
V_\lambda = \bigoplus_{\varpi(\bm \Lambda) = \lambda} V_{\bm \Lambda},
\end{equation*}
where $V_{\bm \Lambda}$ is the subspace of $V_\lambda$ such that for any $v \in V_{\bm \Lambda}$ there is $p \in \bbZ_{> 0}$ such that
\begin{equation*}
(\phi^+_{i, \, m} - \Lambda^+_{i, \, m})^p v = 0, \qquad (\phi^-_{i, \, - m} - \Lambda^-_{i, \, - m})^p v = 0
\end{equation*}
for all $i = 1, \ldots, l$ and $m \in \bbZ_{> 0}$. The space $V_{\bm \Lambda}$ is called the $\ell$-weight space of $\ell$-weight $\bm \Lambda$. We say that $\bm \Lambda$ is an $\ell$-weight of $V$ if $V_{\bm \Lambda} \ne \{0\}$. A nonzero element $v \in V_{\bm \Lambda}$ such that
\begin{equation*}
\phi^+_{i, \, m} v = \Lambda^+_{i, \, m} v, \qquad \phi^-_{i, \, - m} v = \Lambda^-_{i, \, - m} v
\end{equation*}
for all $i = 1, \ldots, l$ and $m \in \bbZ_{> 0}$ is said to be an $\ell$-weight vector of $\ell$-weight $\bm \Lambda$. Every nontrivial $\ell$-weight space contains an $\ell$-weight vector.

For any two $\ell$-weights $\bm \Lambda = (\lambda, \, \bm \Lambda^+ \! , \, \bm \Lambda^-)$ and $\bm \Xi = (\xi, \, \bm \Xi^+ \! , \, \bm \Xi^-)$ we define the $\ell$-weight $\bm \Lambda \, \bm \Xi$ as
\begin{equation}
\bm \Lambda \, \bm \Xi = (\lambda + \xi, \, (\bm \Lambda \, \bm \Xi)^+, \, (\bm \Lambda \, \bm \Xi)^-), \label{lm}
\end{equation}
where
\begin{equation*}
(\bm \Lambda \, \bm \Xi)^+ = (\Lambda^+_i(u) \, \Xi^+_i(u))_{i = 1}^l, \qquad (\bm \Lambda \, \bm \Xi)^- = (\Lambda^-_i(u^{-1}) \, \Xi^-_i(u^{-1}))_{i = 1}^l.
\end{equation*}
The product (\ref{lm}) is an associative operation with respect to which $\gothh^*_\ell$ is an abelian group. Here
\begin{equation*}
\bm \Lambda^{-1} = (- \lambda, (\bm \Lambda^+)^{-1}, \, (\bm \Lambda^-)^{-1}),
\end{equation*}
where
\begin{equation*}
(\bm \Lambda^+)^{-1} = (\Lambda^+_i(u)^{-1})^l_{i = 1}, \qquad (\bm \Lambda^-)^{-1} = (\Lambda^-_i(u)^{-1})^l_{i = 1},
\end{equation*}
and the role of the identity element is played by the $\ell$-weight $(0, \, (\underbracket[.5pt]{1, \, \ldots, \, 1}_l), \, (\underbracket[.5pt]{1, \, \ldots, \, 1}_l))$. Note that by definition each $\Lambda_i^+(u)$ and each $\Lambda_i^-(u^{-1})$ is an invertible power series.

A $\uqlg$-module $V$ in category $\calO$ is called a highest $\ell$-weight module of highest $\ell$-weight $\bm \Lambda$ if there exists an $\ell$-weight vector $v \in V$ of $\ell$-weight $\bm \Lambda$ such that
\begin{equation*}
e_i \, v = 0, \quad i = 1, \ldots, l, \quad \mbox{and} \quad V = \uqlg \, v.
\end{equation*}
The vector with the above properties is unique up to a scalar factor. We call it the  highest $\ell$-weight vector of $V$. Let $V$ and $W$ be highest $\ell$-weight $\uqlg$-modules in category $\calO$ of highest $\ell$-weights $\bm \Lambda$ and $\bm \Xi$ respectively. The submodule of $V \otimes_\Delta W$ generated by the tensor product of the highest $\ell$-weight vectors is a highest $\ell$-weight module of highest $\ell$-weight $\bm \Lambda \, \bm \Xi$.

\subsection{\texorpdfstring{Evaluation representations}{Evaluation representations}} \label{s:lwer}

To construct representations of the quantum loop algebra $\uqlsllpo$, it is common to use the Jimbo's homomorphism $\varepsilon$ from $\uqlsllpo$ to the quantum group $\uqgllpo$ defined by the equations \cite{Jim86a}
\begin{align*}
& \varepsilon(q^{\nu \svee[0]{\alpha}}) = q^{\nu (K_{l+1} - K_1)}, &&  \varepsilon(q^{\nu \svee[i]{\alpha}}) = q^{\nu (K_{i} - K_{i+1})}, \\
& \varepsilon(e_0) = F_{1, \, l + 1} \, q^{K_1 + K_{l+1}}, && \varepsilon(e_i) = E_{i, \, i + 1}, \\
& \varepsilon(f_0) = E_{1, \, l + 1} \, q^{-K_1 - K_{l+1}}, && \varepsilon(f_i) = F_{i, \, i + 1},   
\end{align*}
where $i = 1, \ldots, l$. If $\pi$ is a representation of $\uqgllpo$, then $\pi \circ \varepsilon$ is a representation of $\uqlsllpo$ called an evaluation representation. Given $\mu \in \gothk^*$, we denote the representation $\widetilde \pi^\mu \circ \varepsilon$, where the representation  $\widetilde \pi^\mu$ is described in section \ref{s:vm}, by $\widetilde \varphi^\mu$. Since the representation space of $\widetilde \pi^\mu$ and $\widetilde \varphi^\mu$ is the same, slightly abusing notation, we denote the both corresponding modules as $\widetilde V^\mu$. The same convention is used for the respective finite-dimensional counterparts, if any. It is common to call $\widetilde V^\mu$ an evaluation module. It is a highest weight $\uqlsllpo$-module of highest weight being the restriction of $\mu$ to $\gothh^*$.

In fact, for any $\mu$ in $\gothh^*$, the module $\widetilde V^\mu_\zeta$ is also a highest $\ell$-weight module in category $\calO$. The $\ell$-weight spaces and corresponding $\ell$-weights for $l = 1, 2$ was found in the paper \cite{BooGoeKluNirRaz16}. It appears that all $\ell$-weight spaces are one-dimensional and, therefore, generated by $\ell$-weight vectors.  Although the $\ell$-weight vectors do not coincide with the basis vectors $v^\mu_{\bm m}$, defined by equation (\ref{vml}), they can also be labeled by the $(l + 1)l/2$-tuple $\bm m$, and we use for the $\ell$-weights the notation
\begin{equation*}
\bm \Lambda^\mu_{\bm m}(\zeta) = (\lambda^\mu_{\bm m}, \, \bm \Lambda_{\bm m}^{\mu +}(\zeta), \, \bm \Lambda_{\bm m}^{\mu -}(\zeta)).
\end{equation*}
The highest $\ell$-weight vectors correspond to $\bm m = \bm 0$.

For $l = 1$ we have $\bm m = (m_{1 2})$. The component $\lambda^\mu_{\bm m}$ of the $\ell$-weight $\bm \Lambda^\mu_{\bm m}(\zeta)$ is given by the equation
\begin{equation}
\lambda^\mu_{\bm m} = (\mu_1 - \mu_2 - 2 m_{1 2}) \, \omega_1 \label{psimu}
\end{equation}
and for the only components of $\bm \Lambda_{\bm m}^{\mu +}(\zeta)$ and $\bm \Lambda_{\bm m}^{\mu -}(\zeta)$ we have the expression
\begin{align}
& \Lambda_{\bm m, \, 1}^{\mu +}(\zeta, u) = \frac{1 - q^{2 \mu_1 + 2} \, \zeta^s u}{1 - q^{2 \mu_1 + 2 - 2 m_{1 2}} \, \zeta^s u} \, \frac{1 - q^{2 \mu_2} \, \zeta^s u}{1 - q^{2 \mu_1 - 2 m_{1 2}} \,\zeta^s u},  \label{psimup} \\
& \Lambda_{\bm m, \, 1}^{\mu -}(\zeta, u) = \frac{(1 - q^{- 2 \mu_1 - 2} u^{- 1})(1 - q^{- 2 \mu_2} u^{- 1})}{(1 - q^{- 2 \mu_1 - 2 + 2 m} u^{- 1})(1 - q^{- 2 \mu_1 + 2 m} u^{- 1})}.
\end{align}

For $l = 2$ we have $\bm m = (m_{1 2}, \, m_{1 3}, \, m_{2 3})$. In this case
\begin{equation}
\lambda_{\bm m}^\mu = (\mu_1 - \mu_2 - 2 m_{1 2} - m_{1 3} + m_{2 3}) \, \omega_1 + (\mu_2 - \mu_3 + m_{1 2} - m_{1 3} - 2 m_{2 3}) \, \omega_2, \label{psimu2}
\end{equation}
and for the components of $\bm \Lambda_{\bm m}^{\mu +}(\zeta)$ and $\bm \Lambda_{\bm m}^{\mu -}(\zeta)$ we have
\begin{align}
& \Lambda^{\mu +}_{\bm m, \, 1}(\zeta, u) = \frac{1 - q^{2 \mu_1 - 2 m_{1 3} + 2} \zeta^s u}{1 - q^{2 \mu_1 - 2 m_{1 2} - 2 m_{1 3} + 2} \zeta^s u} \, \frac{1 - q^{2 \mu_2 - 2 m_{2 3}} \zeta^s u}{1 - q^{2 \mu_1 - 2 m_{1 2} - 2 m_{1 3}} \zeta^s u}, \label{psimup2a} \\
& \Lambda^{\mu +}_{\bm m, \, 2}(\zeta, u) = \frac{1 - q^{2 \mu_1 - 2 m_{1 2} - 2 m_{1 3} + 1} \zeta^s u}{1 - q^{2 \mu_1 - 2 m_{1 3} + 1} \zeta^s u} \, \frac{1 - q^{2 \mu_1 + 3} \zeta^s u}{1 - q^{2 \mu_1 - 2 m_{1 3} + 3} \zeta^s u} \notag \\*
& \hspace{13em} {} \times \frac{1 - q^{2 \mu_2 + 1} \zeta^s u}{1 - q^{2 \mu_2 - 2 m_{2 3} + 1} \zeta^s u} \, \frac{1 - q^{2 \mu_3 - 1} \zeta^s u}{1 - q^{2 \mu_2 - 2 m_{2 3} - 1} \zeta^s u}, \label{psimup2b} \\
& \Lambda^{\mu -}_{\bm m, \, 1}(\zeta, u^{-1}) = \frac{1 - q^{- 2 \mu_1 + 2 m_{1 3} - 2} \zeta^s u^{-1}}{1 - q^{- 2 \mu_1 + 2 m_{1 2} + 2 m_{1 3} - 2} \zeta^s u^{-1}} \, \frac{1 - q^{- 2 \mu_2 + 2 m_{2 3}} \zeta^s u^{-1}}{1 - q^{- 2 \mu_1 + 2 m_{1 2} + 2 m_{1 3}} \zeta^s u^{-1}}, \\
& \Lambda^{\mu -}_{\bm m, \, 2}(\zeta, u^{-1}) = \frac{1 - q^{- 2 \mu_1 + 2 m_{1 2} + 2 m_{1 3} - 1} \zeta^s u^{-1}}{1 - q^{- 2 \mu_1 + 2 m_{1 3} - 1} \zeta^s u^{-1}} \, \frac{1 - q^{- 2 \mu_1 - 3} \zeta^s u^{-1}}{1 - q^{- 2 \mu_1 + 2 m_{1 3} - 3} \zeta^s u^{-1}} \notag \\*
& \hspace{13em} {} \times \frac{1 - q^{- 2 \mu_2 - 1} \zeta^ s u^{-1}}{1 - q^{- 2 \mu_2 + 2 m_{2 3} - 1} \zeta^s u^{-1}} \, \frac{1 - q^{- 2 \mu_3 + 1} \zeta^s u^{-1}}{1 - q^{- 2 \mu_2 + 2 m_{2 3} + 1} \zeta^s u^{-1}}.
\end{align}

\subsection{\texorpdfstring{$\ell$-weights of $\uqlbp$-modules}{l-weights of UqL(b+)-modules}}

There are two standard Borel subalgebras of the quantum loop algebra $\uqlg$. In terms of the Drinfeld--Jimbo generators they are defined as follows. The Borel subalgebra $\uqlbp$ is the subalgebra generated by $e^{}_i$ with $0 \le i \le l$ and $q^h$ with $h \in \tgothh$, and the Borel subalgebra $\uqlbm$ is the subalgebra generated by $f^{}_i$ with $0 \le i \le l$ and $q^h$ with $h \in \tgothh$. It is clear that these subalgebras are Hopf subalgebras of $\uqlg$. The description of $\uqlbp$ and $\uqlbm$ in terms of the Drinfeld generators is more intricate.

The category $\calO$ of representations of $\uqlbp$ is defined by exactly the same words as it is defined in section \ref{ss:lw} for the case of $\uqlg$. By definition, any $\uqlbp$-module $V$ in category $\calO$ allows for the weight decomposition
\begin{equation}
V = \bigoplus_{\lambda \in \gothh^*}  V_\lambda, \label{wwp}
\end{equation}
which can be again refined by considering $\ell$-weights.

The Borel subalgebra $\uqlbp$ does not contain the elements $\phi^-_{i, \, - m}$, $i = 1, \ldots, l$, $m \in \bbZ_{>0}$, and its infinite-dimensional commutative subalgebra is generated only by the elements $\phi^+_{i, \, m}$ $i = 1, \ldots, l$, $m \in \bbZ_{>0}$, and $q^h$, $h \in \widetilde \gothh$. Respectively, now an $\ell$-weight $\bm \Lambda$ is a pair
\begin{equation}
\bm \Lambda = (\lambda, \, \bm \Lambda^+), \label{llp}
\end{equation}
where $\lambda \in \gothh^*$ and $\bm \Lambda^+$ is an $l$-tuple $ \bm \Lambda^+ = (\Lambda^+_i(u))_{i = 1}^l$ of formal series
\begin{equation*}
\Lambda^+(u) = 1 + \sum_{m = 1}^\infty \Lambda^+_{i, \, n} \, u^n \in \bbC[[u]].
\end{equation*}
We denote by $\gothh^{* +}_\ell$ the set of $\ell$-weights (\ref{llp}).

For any two $\ell$-weights $\bm \Lambda = (\lambda, \, \bm \Lambda^+)$ and $\bm \Xi = (\xi, \, \bm \Xi^+)$ we define the $\ell$-weight $\bm \Lambda \, \bm \Xi$ as
\begin{equation}
\bm \Lambda \, \bm \Xi = (\lambda + \xi, \, (\bm \Lambda \, \bm \Xi)^+), \label{pk}
\end{equation}
where
\begin{equation*}
(\bm \Lambda \, \bm \Xi)^+ = (\Lambda^+_i(u) \, \Xi^+_i(u))_{i = 1}^l.
\end{equation*}
The product (\ref{pk}) is an associative operation with respect to which $\gothh^{* +}_\ell$ is an abelian group. Here
\begin{equation*}
\bm \Lambda^{-1} = (- \lambda, (\bm \Lambda^+)^{-1}),
\end{equation*}
where
\begin{equation*}
(\bm \Lambda^+)^{-1} = (\Lambda^+_i(u)^{-1})_{i = 1}^l,
\end{equation*}
and the role of the identity element is played by the $\ell$-weight $(0, \, (\underbracket[.5pt]{1, \, \ldots, \, 1}_l))$.

Define a surjective homomorphism $\varpi^+ \colon \gothh^{* +}_\ell \to \gothh^*$ by the relation
\begin{equation*}
\varpi^+(\bm \Lambda) = \lambda
\end{equation*}
if $\bm \Lambda = (\lambda, \, \bm \Lambda^+)$. Now for any $V_\lambda$ entering the decomposition (\ref{wwp}) we have
\begin{equation*}
V_\lambda = \bigoplus_{\varpi^+(\bm \Lambda) = \lambda} V_{\bm \Lambda},
\end{equation*}
where $V_{\bm \Lambda}$ is the subspace of $V_\lambda$ such that for any $v$ in $V_{\bm \Lambda}$ there is $p \in \bbZ_{> 0}$ such that
\begin{equation*}
(\phi^+_{i, \, m} - \Lambda^+_{i, \, m})^p v = 0
\end{equation*}
for all $i = 1, \ldots, l$ and $m \in \bbZ_{> 0}$. The space $V_{\bm \Lambda}$ is called the $\ell$-weight space of $\ell$-weight $\bm \Lambda$. We say that $\bm \Lambda$ is an $\ell$-weight of $V$ if $V_{\bm \Lambda} \ne \{0\}$. A nonzero element $v \in V_{\bm \Lambda}$ such that
\begin{equation*}
\phi^+_{i, \, m} v = \Lambda^+_{i, \, m} v
\end{equation*}
for all $i = 1, \ldots, l$ and $m \in \bbZ_{> 0}$ is said to be an $\ell$-weight vector of $\ell$-weight $\bm \Lambda$. Every nontrivial $\ell$-weight space contains an $\ell$-weight vector.

A $\uqlbp$-module $V$ in is called a highest $\ell$-weight module of highest $\ell$-weight $\bm \Lambda$ if there exists an $\ell$-weight vector $v \in V$ of $\ell$-weight $\bm \Lambda$ such that
\begin{equation*}
e_i \, v = 0, \quad 1 \le i \le l, \quad \mbox{and} \quad V = \uqlbp \, v.
\end{equation*}
The vector with the above properties is unique up to a scalar factor. We call it the  highest $\ell$-weight vector of $V$. Let $V$ and $W$ be highest $\ell$-weight $\uqlbp$-modules in category $\calO$ of highest $\ell$-weights $\bm \Lambda$ and $\bm \Xi$ respectively. The submodule of $V \otimes_\Delta W$ generated by the tensor product of the highest $\ell$-weight vectors is a highest $\ell$-weight module of highest $\ell$-weight $\bm \Lambda \, \bm \Xi$.

For any $\uqlbp$-module $V$ in category $\calO$ and an element $\xi \in \gothh^*$, we define the shifted $\uqlbp$-module $V[\xi]$ shifting the action of the generators $q^h$. Namely, if $\varphi$ is the representation of $\uqlbp$ corresponding to the module $V$ and $\varphi[\xi]$ is the representation corresponding to the module $V[\xi]$, then
\begin{equation}
\varphi[\xi](e_i) = \varphi(e_i), \quad i = 1, \ldots, l, \qquad \varphi[\xi](q^h) = q^{\langle \widetilde \xi, 
\, h \rangle} \varphi(q^h), \quad h \in \tgothh. \label{phixi}
\end{equation}
It is clear that the module $V[\xi]$ is in category $\calO$.

For an arbitrary $\uqlg$-module $V$ we define the shifted module $V[\xi]$, $\xi \in \gothh^*$, as a $\uqlbp$-module obtained first by restricting $V$ to $\uqlbp$ and then shifting the obtained $\uqlbp$-module.

%\begin{additions}
%
%Denoting the $\uqlbp$-module $L(\bm \Lambda_\xi)$ by $[\xi]$, we see that the notations $V[\xi]$ and $\varphi[\xi]$ can be considered as an abbreviation for $V \otimes_\Delta [\xi]$ and $\varphi \otimes_\Delta [\xi]$.
%
%\end{additions}

\subsection{Oscillator representations} \label{s:or}

One can get representations of $\uqlbp$ by restricting to it the representations of $\uqlg$. However, no less important for the theory of integrable systems are representations which cannot be obtained by such a procedure.

In the present paper for the case $\gothg = \sllpo$, we use the family of representations of $\uqlbp$, called the $q$-oscillator representations, which cannot be extended to representations of $\uqlsllpo$. We define them in two steps. First, we introduce a homomorphism of $\uqlbp$ into an $l$th tensor power of the $q$-oscillator algebra $\Osc_q$, then for each factor of the tensor product we choose a suitable representation of $\Osc_q$.

The $q$-oscillator algebra $\Osc_q$ is an algebra with generators $b^\dagger$, $b$ and $q^{\nu N}$, $\nu \in \bbC$, and relations
\begin{gather*}
q^0 = 1, \qquad q^{\nu_1 N} q^{\nu_2 N} = q^{(\nu_1 + \nu_2)N}, \\
q^{\nu N} b^\dagger q^{-\nu N} = q^\nu b^\dagger, \qquad q^{\nu N} b q^{-\nu N} = q^{-\nu} b, \\
b^\dagger b = \frac{q^N - q^{-N}}{q - q^{-1}}, \qquad b b^\dagger = \frac{q q^N - q^{-1} q^{-N}}{q - q^{-1}}.
\end{gather*}
Two basic representations of $\Osc_q$ are of particular interest. First, let $W^+$ be the free vector space on the sequence $(w_n)_{n \in \bbZ_{ \ge 0}}$. One can show that the relations
\begin{gather*}
q^{\nu N} w_n = q^{\nu n} w_n,  \\[.3em]
b^\dagger w_n = w_{n + 1}, \qquad b \, w_n = [n]_q w_{n - 1},
\end{gather*}
where we assume that $w_{-1} = 0$, endow $W^+$ with the structure of an $\Osc_q$-module. We denote the corresponding representation of the algebra $\Osc_q$ by $\chi^+$. Further, let $W^-$ be again a free vector space on the sequence $(w_n)_{n \in \bbZ_{ \ge 0}}$. The relations
\begin{gather*}
q^{\nu N} w_n = q^{- \nu (n + 1)} v_n,  \\[.3em]
b \, w_n = w_{n + 1}, \qquad b^\dagger w_n = - [n]_q w_{n - 1},
\end{gather*}
where we again assume that $w_{-1} = 0$, endow $W^-$ with the structure of an $\Osc_q$-module. We denote the corresponding representation of $\Osc_q$ by $\chi^-$.

We consider the tensor product of $l$ copies of the $q$-oscillator algebra, and denote
\begin{gather*}
b_i = \underbracket[.5pt]{1 \otimes \cdots \otimes 1}_{i - 1} {} \otimes b \otimes \underbracket[.5pt]{1 \otimes \cdots \otimes 1}_{l - i} \, , \qquad
b_i^\dagger =  \underbracket[.5pt]{1 \otimes \cdots \otimes 1}_{i - 1} {} \otimes b^\dagger \otimes \underbracket[.5pt]{1 \otimes \cdots \otimes 1}_{l - i} \, , \\
q^{\nu N_i} =  \underbracket[.5pt]{1 \otimes \cdots \otimes 1}_{i - 1} {} \otimes q^{\nu N} \otimes \underbracket[.5pt]{1 \otimes \cdots \otimes 1}_{l - i} \, .
\end{gather*}
In the paper \cite{NirRaz17b} the homomorphism $o$ of $\uqlbp$ into $(\Osc_q)^{\otimes \, l}$ 
described by the equations
\begin{align*}
& o(q^{\nu \svee[0]{\alpha}}) = q^{\nu (2 N_1 + \sum_{j = 2}^l N_j)}, 
&& o(e_0) = b^\dagger_1 \, q^{\sum_{j = 2}^l N_j}, \\[.3em]
& o(q^{\nu \svee[i]{\alpha}}) = q^{\nu (N_{i + 1} - N_{i})}, 
&& o(e_i) = - b^{}_i \, b^{\mathstrut \dagger}_{i + 1} \, q^{N_i - N_{i + 1} - 1}, \\[.3em]
& o(q^{\nu \svee[l]{\alpha}}) = q^{- \nu (2 N_l + \sum_{j = 1}^{l - 1} N_j)},
&& o(e_l) = - (q - q^{-1})^{-1} \, b^{}_l \, q^{N_l},
\end{align*}
where $1 \le i < l$, was obtained by some limiting procedure starting from the evaluation representations $\widetilde \varphi^\lambda$ of $\uqlsllpo$. We define
\begin{equation*}
\theta = (\underbracket[.5pt]{\, \chi^+ \otimes \cdots \otimes \chi^+}_l) \circ o.
\end{equation*}
The basis of this representation is formed by the vectors
\begin{equation*}
w_{\bm n} = b^{\dagger \, n_1}_1 \ldots b^{\dagger \, n_l}_l \, w^{}_{\bm 0},
\end{equation*}
where $n_i \in \bbZ_{\ge 0}$ for all $1 \le i \le l$, and we use the notation $\bm n = (n_1, \ldots, n_l)$. Here the vector $w_{\bm 0}$ is the vacuum vector, satisfying the equations
\begin{equation*}
b_i^{} \, w_{\bm 0} = 0, \qquad 1 \le i \le l.
\end{equation*}
The representation $\theta$ is in category $\calO$. It is a highest $\ell$-weight representation with the highest $\ell$-weight vector $w_{\bm 0}$.

There is an automorphism of $\uqlsllpo$ defined by the equations
\begin{equation*}
\sigma(q^{\nu \svee[i]{\alpha}}) = q^{\nu \svee[i + 1]{\alpha}}, \qquad \sigma(e_i) = e_{i + 1}, \qquad \sigma(f_i) = f_{i + 1}, \qquad 0 \le i \le l,
\end{equation*}
where it is assumed that $q^{\nu \svee[l + 1]{\alpha}} = q^{\nu \svee[0]{\alpha}}$, $e_{l + 1} = e_0$ and $f_{l + 1} = f_0$. One can restrict $\sigma$ to an automorphism of $\uqlbp$. It is useful to have in mind that $\sigma^{l + 1} = 1$. We define a collection of homomorphisms from $\uqlbp$ into $\Osc_q^{\otimes l}$ as
\begin{equation*}
o_a = o \circ \sigma^{-a},
\end{equation*}
and a family of representations
\begin{equation}
\theta_a = \chi_a \circ o_a, \qquad a = 1, \ldots, l + 1,
\label{thetaa}
\end{equation}
where
\begin{equation*}
\chi_a = \underbracket[.5pt]{\, \chi^- \otimes \cdots \otimes \chi^-}_{l - a + 1} \otimes \underbracket[.5pt]{\, \chi^+ \otimes \cdots \otimes \chi^+}_{a - 1}.
\end{equation*}
Note that $\theta = \theta_{l + 1}$. The corresponding basis vectors are 
\begin{equation*}
w_{a, \, \bm n} = b_1^{n_1} \ldots b_{l - a + 1}^{n_{l - a + 1}} \, b^{\dagger \, n_{l - a + 2}}_{l - a + 2} \ldots b^{\dagger \, n_l}_{l} \, w_{a, \, \bm 0}.
\end{equation*}
The vacuum vector $w_{a, \, \bm 0}$ satisfies the equations
\begin{equation*}
b_i^\dagger w_{a, \, \bm 0} = 0, \quad 1 \le i \le l - a + 1, \qquad b_i^{} \,w_{a, \, \bm 0} = 0, \quad l - a + 2 \le i \le l.
\end{equation*}
All the representations $\theta_a$ are highest $\ell$-weight representations in category $\calO$. The vectors $w_{a, \, \bm n}$ are the $\ell$-weight vectors, and $w_{a, \, \bm 0}$ is the highest $\ell$-weight vector of the representation $\theta_a$.

The explicit form of the $\ell$-weights $\bm \Psi_{a, \, \bm n}(\zeta) = (\psi_{a, \, \bm n}, \, \bm \Psi^+_{a, \, \bm n}(\zeta))$ for the representations $(\theta_a)_\zeta$ was found in the paper~\cite{BooGoeKluNirRaz17b}. For future usage, we give below the relevant expressions for $l = 1$ and $l = 2$. For $l = 1$ we have $\bm n = (n_1)$ and
\begin{align*}
& \psi_{1, \, \bm n} = - (2 n_1 + 2) \omega_1, \\
& \Psi_{1, \, \bm n, \, 1}(\zeta, \, u) = \frac{1 - q \zeta^s u}{(1 - q^{-2 n_1 - 1} \zeta^s u) (1 - q^{-2 n_1 + 1} \zeta^s u)}, \\
& \psi_{2, \, \bm n} = - 2 n_1 \omega_1, \\
& \Psi_{2, \, \bm n, \, 1}(\zeta, \, u) = 1 - q \zeta^s u.
%\end{align*}
\intertext{For $l = 2$ we have $\bm n = (n_1, \, n_2)$ and}
%\begin{align*}
& \psi_{1, \, \bm n} = (-2 n_1 - n_2 - 3) \omega_1 + (n_1 - n_2) \omega_2, \\
& \Psi_{1, \, \bm n, \, 1}(\zeta, \, u) = \frac{1 - q^{- 2 n_2} \zeta^s u}{(1 - q^{-2 n_1 - 2 n_2} \zeta^s u) (1 - q^{-2 n_1 - 2 n_2 - 2} \zeta^s u)}, \\
& \Psi_{1, \, \bm n, \, 2}(\zeta, \, u) = \frac{(1 - q \zeta^s u) (1 - q^{- 2 n_1 -2 n_2 - 1} \zeta^s u)}{(1 - q^{-2 n_2 - 1} \zeta^s u) (1 - q^{- 2 n_2 + 1} \zeta^s u)}, \\
& \psi_{2, \, \bm n} = (n_1 - 2 n_2 + 1) \omega_1 + (- 2 n_1 + n_2 - 2) \omega_2, \\
& \Psi_{2, \, \bm n, \, 1}(\zeta, \, u) = 1 - q^{- 2 n_1} \zeta^s u, \\
& \psi_{2, \, \bm n, \, 1}(\zeta, \, u) = \frac{1 - q \zeta^s u}{(1 - q^{- 2 n_1 - 1} \zeta^s u)(1 - q^{- 2 n_1 +1})}, \\
& \psi_{3, \, \bm n} = (- n_1 + n_2) \omega_1 + (- n_1 - 2 n_2) \omega_2, \\
& \Psi_{3, \, \bm n, \, 1}(\zeta, \, u) = 1, \\
& \Psi_{3, \, \bm n, \, 2}(\zeta, \, u) = 1 - q \zeta^s u,
\end{align*}

\subsection{\texorpdfstring{$q$-characters and Grothendieck ring}{q-characters and Grothendieck ring}}

Let $V$ be a $\uqlg$-module in category $\calO$. We define the character of $V$ as a formal sum
\begin{equation*}
\ch(V) = \sum_{\lambda \in \gothh^*} \dim \, V_\lambda \, [\lambda].
\end{equation*}
By the definition of the category $\calO$,  $\dim V_\lambda = 0$ for $\lambda$ outside the union of a finite number of sets of the form $D(\mu)$, $\mu \in \gothh^*$. For any two $\uqlg$-modules $V$ and $U$ in category $\calO$ we have
\begin{equation*}
\ch(V \oplus U) = \ch(V) + \ch(U).
\end{equation*}
More generally, if $\uqlg$-modules $V$, $W$ and $U$ in category $\calO$ can be included in a short exact sequence
\begin{equation}
0 \rightarrow V \rightarrow W \rightarrow U \rightarrow 0, \label{vwu}
\end{equation}
then
\begin{equation*}
\ch(W) = \ch(V) + \ch(U).
\end{equation*}
It can be also shown that
\begin{equation*}
\ch(V \otimes_\Delta U) = \ch(V) \, \ch(U)
\end{equation*}
for any $\uqlg$-modules $V$ and $U$ in category $\calO$. Here, to multiply characters we assume that
\begin{equation*}
[\lambda] \, [\mu] = [\lambda + \mu]
\end{equation*}
for any $\lambda, \mu \in \gothh$.

The Grothendieck group of the category $\calO$ of $\uqlg$-modules is defined as the quotient of the free abelian group on the set of all isomorphism classes of objects in $\calO$ modulo the relations
\begin{equation*}
\langle V \rangle = \langle U \rangle + \langle W \rangle
\end{equation*}
if the objects $V$, $W$ and $U$ can be included in a short exact sequence (\ref{vwu}). Here for any object $V$, $\langle V \rangle$ denotes the isomorphism class of $V$. Defining 
\begin{equation*}
\langle V \rangle \langle W \rangle = \langle V \otimes_\Delta W \rangle,
\end{equation*}
we come to the Grothendieck ring of $\calO$. It is a commutative unital ring, for which the role of the unit is played by the trivial $\uqlg$-module. We see that the character can be considered as a mapping from the Grothendieck ring of $\calO$. However, it is not injective, and, therefore, does not uniquely distinguish between its elements. This role is played by the $q$-character.

We define the $q$-character of a $\uqlg$-module $V$ in category $\calO$ as
\begin{equation*}
\ch_q(V) = \sum_{\bm \Lambda \in \gothh^*_\ell} \dim (V_{\bm \Lambda}) \, [\bm \Lambda].
\end{equation*}
One can easily demonstrate that
\begin{equation*}
\varpi(\ch_q(V)) = \ch(V).
\end{equation*}
Here we assume that
\begin{equation*}
\varpi([\bm \Lambda]) = [\varpi(\bm \Lambda)],
\end{equation*}
and extend this rule by linearity. The $q$-character has the same properties as the usual character. Namely, if $\uqlg$-modules $V$, $W$ and $U$ in category $\calO$ can be included in a short exact sequence (\ref{vwu}), then
\begin{equation*}
\ch_q(W) = \ch_q(V) + \ch_q(U),
\end{equation*}
and for any two $\uqlg$-modules $V$ and $U$ in category $\calO$, we have
\begin{equation}
\ch_q(V \otimes_\Delta W) = \ch_q(V) \, \ch_q(W), \label{cvtu}
\end{equation}
see the paper \cite{FreRes99}. To define the product of $q$-characters, we assume that
\begin{equation*}
[\bm \Lambda][\bm \Xi] = [\bm \Lambda \bm \Xi]
\end{equation*}
for any $\bm \Lambda, \bm \Xi \in \gothh^*_\ell$. Thus, $q$-character can also be considered as a mapping from the Grothendieck ring of $\calO$, and now one can demonstrate that this mapping is injective. In other words, different elements of the Grothendieck ring of $\calO$ have different $q$-characters.

It follows from (\ref{cvtu}) that
\begin{equation*}
\ch_q(V \otimes_\Delta W) = \ch_q(V) \ch_q(W) = \ch_q(W) \ch_q(V) = \ch_q(W \otimes_\Delta V).
\end{equation*}
It means that the $\uqlg$-modules $V \otimes_\Delta W$ and $W \otimes_\Delta V$ belong to the same equivalence class in the Grothendieck ring of $\calO$.

All above can be naturally extended to the categories $\calO$ of the modules over $\uqlbp$ and $\uqlbm$.

\section{\texorpdfstring{Universal $R$-matrix and integrability objects}{Universal R-matrix and integrability objects}} \label{s:urm}

%\subsection{\texorpdfstring{Universal $R$-matrix}{Universal R-matrix}}

\subsection{\texorpdfstring{Quantum group as a $\bbC[[\hbar]]$-algebra}{Quantum group as a C[[h]]-algebra}}

As any Hopf algebra a quantum loop algebra $\uqlg$ has another comultiplication called the opposite comultiplication. It is defined by the equation
\begin{equation*}
\Delta' = \Pi \circ \Delta,
\end{equation*}
where
\begin{equation*}
\Pi(x \otimes y) = y \otimes x
\end{equation*}
for all $x, y \in \uqlg$.

There are several different definitions of quantum groups. In particular, $\hbar$ can be not only a complex number \cite{Jim86a, EtiFreKir98, Dam98}, but also an indeterminate, so that the quantum group is a $\bbC[[\hbar]]$-algebra \cite{Dri87, TolKho92, KhoTol92, KhoTol93, KhoTol94}. Let us assume temporally that it is the case. Herewith $\uqlg$ is a quasitriangular Hopf algebra. It means that, up to a central element, there exists a unique invertible element $\calR$ of the completed tensor product $\uqlg \mathop{\widehat \otimes} \uqlg$, called the universal $R$-matrix, such that
\begin{equation}
\Delta'(x) = \calR \, \Delta(x) \, \calR^{-1} \label{urmd}
\end{equation}
for all $x \in \uqlg$, and\footnote{For an explanation of the notation see, for example, the paper \cite{NirRaz19}.}
\begin{equation}
(\Delta \otimes \id) (\calR) = \calR^{(1 3)} \calR^{(2 3)}, \qquad (\id \otimes \Delta) (\calR) = \calR^{(1 3)} \calR^{(1 2)}. \label{urm}
\end{equation}
The latter equations are considered as equalities in the completed tensor product of three copies of $\uqlg$. In fact, it follows from the explicit expression for the universal $R$-matrix \cite{TolKho92, KhoTol92, KhoTol93, KhoTol94} that it is an element of a completed tensor product of two Borel subalgebras $\uqlbp$ and $\uqlbm$.

A general integrability object is defined as follows. Let $\varphi$ be a representation of the Borel subalgebra $\uqlbp$ on a vector space $V$, and $\psi$ a representation of $\uqlbm$ on a vector space $U$.\footnote{In the present paper, we always assume that the representations used are in category $\calO$.}  The corresponding integrability object $\bm X_{\varphi | \psi}$ is defined by the equation
\begin{equation*}
\rho_{\varphi | \psi} \bm X_{\varphi | \psi} = [(\varphi \otimes \psi)(\calR)],
\end{equation*}
where $\rho_{\varphi | \psi}$ is a scalar normalization factor. It is evident that this object is an element of $\End(V) \otimes \End(U)$. It follows from (\ref{urmd}) that
\begin{equation}
(\varphi \otimes \psi)(\Pi(\Delta(x))) = \bm X_{\varphi | \psi} (\varphi \otimes \psi)(\Delta(x)) (\bm X_{\varphi | \psi})^{-1}, \label{urmx}
\end{equation}
and equation (\ref{urm}) gives
\begin{equation}
\bm X\strut_{\varphi_1 \otimes_\Delta \varphi_2 | \psi} = \bm X_{\varphi_1 | \psi}^{(1 3)}\bm  X_{\varphi_2 | \psi}^{(2 3)}, \qquad \bm X\strut_{\varphi | \psi_1 \otimes_\Delta \psi_2} = \bm X_{\varphi | \psi_2}^{(1 3)} \bm X_{\varphi | \psi_1}^{(1 2)}. \label{xfff}
\end{equation}
In fact, hereinafter, we assume that
\begin{equation*}
\rho_{\varphi_1 \otimes_\Delta \varphi_2 | \psi} = \rho_{\varphi_1 | \psi} \, \rho_{\varphi_2 | \psi}, \qquad \rho_{\varphi_1 | \psi_1 \otimes_\Delta \psi_2} = \rho_{\varphi | \psi_1} \, \rho_{\varphi | \psi_2}.
\end{equation*}

\subsection{\texorpdfstring{Quantum group as a $\bbC$-algebra}{Quantum group as a C-algebra}}

The expression for the universal $R$-matrix of a quantum loop algebra $\uqlg$ considered as a $\bbC[[\hbar]]$-algebra can be constructed using the procedure proposed by Khoroshkin and Tolstoy \cite{TolKho92, KhoTol93, KhoTol94}. However, in this paper, we define the quantum loop algebra $\uqlg$ as a $\bbC$-algebra. In fact, one can use the expression for the universal $R$-matrix from the papers \cite{TolKho92, KhoTol93, KhoTol94} to construct the integrability objects also for this case, having in mind that $\uqlg$ is quasitriangular only in some restricted sense, see the paper \cite{Tan92}, the book \cite[p.~327]{ChaPre94}, and the discussion below.

We restrict ourselves to the following case. Let $\varphi$ be a representation of $\uqlbp$ on a vector space $V$, and $\psi$ a representation of $\uqlbm$ on a vector space $U$. Define an integrability object $\bm X_{\varphi | \psi}$ as an element of $\End(V) \otimes \End(U)$ by the equation
\begin{equation}
\rho_{\varphi | \psi} \bm X_{\varphi | \psi} = (\varphi \otimes \psi)(\calR_{\prec \delta} \, \calR_{\sim \delta} \, \calR_{\succ \delta}) \bm K_{\varphi | \psi}. \label{xfpb}
\end{equation}
Here $\calR_{\prec \delta}$, $\calR_{\sim \delta}$ and $\calR_{\succ \delta}$ are elements of $\uqlbp \widehat \otimes \uqlbm$, $\bm K_{\varphi | \psi}$ is an element of $\End(V) \otimes \End(U)$, and $\rho_{\varphi | \psi}$ a scalar normalization factor.

The explicit expressions for the elements $\calR_{\prec \delta}$ , $\calR_{\sim \delta} $ and $\calR_{\succ \delta}$ for the case $\gothg = \sllpo$ are given in the paper \cite{NirRaz19}. For a general case one can consult the papers \cite{TolKho92, KhoTol93, KhoTol94}. The element $\bm K_{\varphi | \psi}$ for general $\gothg$ is given by the equation
\begin{equation}
\bm K_{\varphi | \psi} = \sum_{\lambda \in \gothh^*} \varphi \big( q^{- \sum_{i, \, j = 1}^l \svee[\raisebox{-.17em}{$\scriptscriptstyle i$}]{\alpha} c^{}_{i j} d_j^{-1} \, \langle \lambda, \, \svee[\raisebox{-.17em}{$\scriptscriptstyle j$}]{\alpha} \rangle} \big) \otimes \Pi_\lambda,
\label{kfpb}
\end{equation}
where $\Pi_\lambda \in \End(U)$ is the projector on the component $U_\lambda$ of the weight decomposition
\begin{equation*}
U = \bigoplus_{\lambda \in \gothh^*}  U_\lambda,
\end{equation*}
and $c_{i j}$ are the matrix entries of the matrix $C$ inverse to the Cartan matrix $A$ of the Lie algebra $\gothg$. Using formulas of the paper \cite{Usm94}, for $\gothg = \gllpo$ we obtain
\begin{equation*}
c_{i j} = \frac{i (l - j + 1)}{l + 1}, \quad i \le j, \qquad c_{i j} = \frac{(l - i + 1) j}{l + 1}, \quad i > j.
\end{equation*}

One can show that the integrability objects defined by equation (\ref{xfpb}) satisfy equations~(\ref{urmx}) and (\ref{xfff}). So the integrability objects defined by equation (\ref{xfpb}) behave as if they were constructed from the real universal $R$-matrix.

Each integrability object $\bm X_{\varphi | \psi} \in \End(V) \otimes \End(U)$ generates an integrability object $\bm Y_{\varphi | \psi} \in \End(U)$ defined as
\begin{equation}
\bm Y_{\varphi | \psi} = (\tr_\varphi \mathop{\otimes} \id_{\End(U)})(\bm X_{\varphi | \psi} (\varphi(q^t) \otimes 1_{\End(U)})), \label{yfp}
\end{equation}
where $q^t \in \uqlg$ is a group like\footnote{An element $x \in \uqlg$ is called group like if $\Delta(x) = x \otimes x$.} twisting element necessary for the convergence of the trace, and $\tr_\varphi$ is the trace on $\uqlbp$ defined by the equation
\begin{equation}
\tr_\varphi = \tr_{\End(V)} \circ \varphi, \label{traa}.
\end{equation}
Here $\tr_{\End(V)}$ is the usual trace on the endomorphism algebra of the vector space $V$. By definition, the integrability object $\bm Y_{\varphi | \psi}$ depends only on the equivalence class in the Grothendieck ring to which the representation $\varphi$ belongs.

It is productive to define the corresponding universal integrability objects. They are defined as formal objects with specific rules of use. If $\varphi$ is a representation of $\uqlbp$ on a vector space $V$, the universal integrability object $\calX_{\varphi}$, corresponding to the integrability objects of type $X$, behaves as an element of $\End(V) \otimes \uqlbm$ and obeys the rule
\begin{equation*}
(\id_{\End(V)} \mathop{\otimes} \psi) (\calX_\varphi) = \rho_{\varphi | \psi} \bm X_{\varphi | \psi}
\end{equation*}
for any representation $\psi$ of $\uqlbm$. The universal integrability object $\calY_{\varphi}$, corresponding to the integrability objects of type $Y$, behaves as an element of $\uqlbm$ and obeys the rule
\begin{equation*}
\psi (\calY_\varphi) = \rho_{\varphi | \psi} \bm Y_{\varphi | \psi}.
\end{equation*}
It follows from (\ref{yfp}) that
\begin{equation*}
\calY_\varphi^\chi = (\tr_\varphi \mathop{\otimes} \id_{\uqlg})(\calX_\varphi (\varphi(q^t) \otimes 1_{\uqlg})).
\end{equation*}

Let $\varphi_1$ and $\varphi_2$ be representations of $\uqlbp$ on vector spaces $V_1$ and $V_2$, and $\psi$ a representation of $\uqlbm$ on a vector space $U$. Define the tensor product of two universal integrability  objects $\calX_{\varphi_1}$ and $\calX_{\varphi_2}$, behaving as an element of $\End(V_1) \otimes \End(V_2) \otimes \uqlbm$ by the rule
\begin{equation*}
(\id_{\End(V_1)} \otimes \id_{\End(V_2)} \otimes \psi)(\calX_{\varphi_1} \otimes \calX_{\varphi_2}) = \rho_{\varphi_1 | \psi} \rho_{\varphi_2 | \psi} (\bm X_{\varphi_1 | \psi})^{(1 2)} (\bm X_{\varphi_2 | \psi})^{(1 3)},
\end{equation*}
and the product of two universal integrability  objects $\calY_{\varphi_1}$ and $\calY_{\varphi_2}$ behaving as an element of $\uqlbm$, by the rule
\begin{equation*}
\psi(\calY_{\varphi_1} \calY_{\varphi_2}) = \rho_{\varphi_1 | \psi} \, \rho_{\varphi_2 | \psi} \bm Y_{\varphi_1 | \psi} \bm Y_{\varphi_2 | \psi}.
\end{equation*}
Using the first equation of (\ref{xfff}), we obtain
\begin{equation*}
\bm Y_{\varphi_1 \otimes_\Delta \varphi_2 | \psi} = \bm Y_{\varphi_1 | \psi} \bm Y_{\varphi_2 | \psi},
\end{equation*}
or, in terms of universal integrability objects,
\begin{equation}
\calY_{\varphi_1 \otimes_\Delta \varphi_2 | \psi} = \calY_{\varphi_1 | \psi} \calY_{\varphi_2 | \psi}.  \label{yffp}
\end{equation}
Since the representations $\varphi_1 \otimes_\Delta \varphi_2$ and $\varphi_2 \otimes_\Delta \varphi_1$ belongs to the same equivalence class in the Grothendieck ring of $\calO$, we have
\begin{equation*}
\bm Y_{\varphi_1 \otimes_\Delta \varphi_2 | \psi} = \bm Y_{\varphi_2 \otimes_\Delta \varphi_1 | \psi},
\end{equation*}
and, therefore,
\begin{equation}
\bm Y_{\varphi_1 | \psi} \bm Y_{\varphi_2 | \psi} = \bm Y_{\varphi_2 | \psi} \bm Y_{\varphi_1 | \psi}, \label{yyyy}
\end{equation}
or, in terms of universal integrability objects,
\begin{equation*}
\calY_{\varphi_1} \calY_{\varphi_2} = \calY_{\varphi_2} \calY_{\varphi_1}. \label{yyyyu}
\end{equation*}

Let now $x \in \uqlbp \cup \uqlbm$ be a group like element commuting with the twisting element $q^t$. Starting with equation (\ref{urmx}), we obtain
\begin{equation*}
\bm Y_{\varphi | \psi} \psi(x) = \psi(x) \, \bm Y_{\varphi | \psi}.
\end{equation*}
In particular, for any $h \in \widetilde \gothh$, we have
\begin{equation*}
\bm Y_{\varphi | \psi} \psi(q^h) = \psi(q^h) \, \bm Y_{\varphi | \psi}.
\end{equation*}
Hence, in terms of universal integrability objects
\begin{equation*}
\calY_\varphi x = x \, \calY_\varphi,
\end{equation*}
and
\begin{equation*}
\calY_\varphi q^h = q^h \, \calY_\varphi.
\end{equation*}

Consider the behavior of an integrability objects $\bm X_{\varphi | \psi}$ and $\calX_{\varphi | \psi}$ under a shift of the homomorphism $\varphi$, see equation (\ref{phixi}). First, since the elements $\calR_{\prec \delta}$, $\calR_{\sim \delta}$ and $\calR_{\succ \delta}$ do not depend on the generators $q^h$, $h \in \widetilde \gothh$, see the papers \cite{TolKho92, KhoTol93, KhoTol94}, we determine that
\begin{equation*}
(\varphi[\xi] \otimes \psi)(\calR_{\prec \delta} \, \calR_{\sim \delta} \, \calR_{\succ \delta}) = (\varphi \otimes \psi)(\calR_{\prec \delta} \, \calR_{\sim \delta} \, \calR_{\succ \delta}).
\end{equation*}
Then, using equations (\ref{phixi}) and (\ref{kfpb}), we obtain
\begin{equation*}
K_{\varphi[\xi] | \psi} = K_{\varphi | \psi} \big( 1_{\End(V)} \otimes \psi(q^{- \sum_{i, \, j = 1}^l \langle \xi, \, \svee[\raisebox{-.17em}{$\scriptscriptstyle i$}]{\alpha} \rangle c^{}_{i j} d^{-1}_j \, \svee[\raisebox{-.17em}{$\scriptscriptstyle j$}]{\alpha}}) \big)
\end{equation*}
and come to the equations
\begin{equation*}
\bm X_{\varphi[\xi] | \psi} = X_{\varphi | \psi} \, \big( 1_{\End(V)} \otimes \psi(q^{- \sum_{i, \, j = 1}^l \langle \xi, \, \svee[i]{\alpha} \rangle c_{i j} \, \svee[j]{\alpha}}) \big), \quad  \calX_{\varphi[\xi]} = \calX_\varphi \, \big( 1_{\End(V)} \otimes q^{- \sum_{i, \, j = 1}^l \langle \xi, \, \svee[i]{\alpha} \rangle c_{i j} \, \svee[j]{\alpha}} \big).
\end{equation*}
For the integrability objects $\bm Y_{\varphi | \psi}$ and $\calY_\varphi$ we have the equations
\begin{equation}
\bm Y_{\varphi[\xi] | \psi} = Y_{\varphi | \psi} \, \psi \big( q^{- \sum_{i, \, j = 1}^l \langle \xi, \, \svee[i]{\alpha} \rangle c_{i j} \, \svee[j]{\tilde \alpha}} \big), \qquad \calY_{\varphi[\xi]} = \calY_\varphi \, q^{-  \sum_{i, \, j = 1}^l \langle \xi, \, \svee[i]{\alpha} \rangle c_{i j} \, \svee[j]{\tilde \alpha}}, \label{yfsp}
\end{equation}
where
\begin{equation*}
\svee[i]{\tilde \alpha} = \svee[i]{\alpha} - t_i.
\end{equation*}
Here and below we set
\begin{equation}
t = \sum_{i = 1}^l t_i \, \svee[i]{\alpha}, \label{tqfh}
\end{equation}
where $t_i$, $i = 1, \ldots, l$ are complex twisting parameters.

In fact, the integrability objects depend on spectral parameters. To define such objects, we use as a representation $\varphi$ the representation $\varphi_\zeta$, and as $\psi$ the $n$th tensor power of a homomorphism $\psi$.\footnote{One can also use as $\psi$ the tensor product $\psi_{\eta_1} \otimes_\Delta \cdots \otimes_\Delta \psi_{\eta_n}$. However, we do not consider such a generalization in this paper.} The corresponding universal integrability objects are denoted as
\begin{equation*}
\calX_\varphi(\zeta) =  \calX_{\varphi_\zeta}, \qquad \calY_\varphi(\zeta) =  \calY_{\varphi_\zeta}
\end{equation*}
while for usual integrability objects we use the notation 
\begin{equation}
\underset{n}{\bm X}{}_{\varphi | \psi} (\zeta) 
= \bm X_{\varphi_\zeta | \underbracket[.5pt]{\scriptstyle \psi \otimes_\Delta \cdots \otimes_\Delta \psi}_n}, \qquad \underset{n}{\bm Y}{}_{\varphi | \psi} (\zeta) 
= \bm Y_{\varphi_\zeta | \underbracket[.5pt]{\scriptstyle \psi \otimes_\Delta \cdots \otimes_\Delta \psi}_n}. \label{xfpz}
\end{equation}
When $n = 1$ we usually write just $\bm X_{\varphi | \psi} (\zeta)$ and $\bm Y_{\varphi | \psi} (\zeta)$. In accordance with our conventions, we set
\begin{equation*}
\underset{n}{\rho}{}_{\varphi | \psi}(\zeta) = \rho_{\varphi_\zeta | \underbracket[.5pt]{\scriptstyle \psi \otimes_\Delta \cdots \otimes_\Delta \psi}_n} = (\rho_{\varphi_{\zeta} | \psi})^n = \rho_{\varphi | \psi}(\zeta)^n.
\end{equation*}

All integrability objects that we use in this paper, are constructed as described above. However, depending on the role they play in the integration procedure, they are given different names. Below we describe the main classes of integrability objects. It is worth to have in mind that the proposed nomination is rather conditional, although it is widespread.

The most famous integrability objects are the $R$ -operators. They form a special class of integrability objects of type $X$ used to permute integrability objects of type $Y$. However, there is a more general method for demonstrating commutativity of integrability objects of type $Y$, described in the previous section, and we do not define and use $R$-operators in the present paper.

\subsection{Monodromy operators and transfer operators}

When $\varphi$ is a representation of the quantum loop algebra $\uqlg$ on a vector space $V$ and $\psi$ a representation of $\uqlg$ on a vector space $U$, the integrability object $\underset{n}{\bm X}{}_{\varphi | \psi}(\zeta)$ is called a monodromy operator and denoted $\underset{n}{\bm M}{}_{\varphi | \psi}(\zeta)$.

The type $Y$ companion of a monodromy operator $\underset{n}{\bm M}{}_{\varphi | \psi}(\zeta)$ is called the transfer operator and denoted $\underset{n}{\bm T}{}_{\varphi | \psi}(\zeta)$. Explicitly, we have
\begin{equation*}
\underset{n}{\bm T}{}_{\varphi | \psi}(\zeta) = \big( \tr_{\End(V)} \otimes \id_{\End(W^{\otimes n})} \big) \big( \underset{n}{\bm M}{}_{\varphi | \psi}(\zeta) (\varphi_\zeta(t) \otimes 1_{\End(W^{\otimes n})}) \big),
\end{equation*}
where $t$ is a twisting element, which we define by equation (\ref{tqfh}).

Let $\varphi_1$ and $\varphi_2$ be representations of $\uqlg$, respectively, and $\psi$ a representation of $\uqlg$. It follows from equation (\ref{yyyy}) that
\begin{equation*}
\underset{n}{\bm T}{}_{\varphi_1 | \psi}(\zeta_1) \, \underset{n}{\bm T}{}_{\varphi_2 | \psi}(\zeta_2) =\underset{n}{\bm T}{}_{\varphi_2 | \psi}(\zeta_2) \, \underset{n}{\bm T}{}_{\varphi_1 | \psi}(\zeta_1),
\end{equation*}
for any $\zeta_1, \zeta_2 \in \bbC^\times$. Similar commutativity takes place for the universal transfer operators,
\begin{equation*}
\calT_{\varphi_1}(\zeta_1) \calT_{\varphi_2}(\zeta_2) = \calT_{\varphi_2}(\zeta_2) \calT_{\varphi_2}(\zeta_2).
\end{equation*}
In fact, it is important for commutativity that the twist element is group like.

In this paper we construct the monodromy operators using as $\varphi$ the evaluation representations $\widetilde \varphi^\lambda$ and $\varphi^\lambda$ defined in section \ref{s:lwer}, and as $\psi$ the $(l + 1)$-dimensional evaluation representation $\varphi^{\omega_1}$. Here the following notation is used
\begin{equation*}
\underset{n}{\widetilde{\bm T}}{}^\lambda(\zeta) = \underset{n}{\widetilde{\bm T}}{}_{\widetilde \varphi^\lambda | \varphi^{\omega_1}}(\zeta), \qquad \underset{n}{\bm T}{}^\lambda(\zeta) = \underset{n}{\bm T}{}_{\varphi^\lambda | \varphi^{\omega_1}}(\zeta).
\end{equation*}
Similarly, for the corresponding universal transfer operators we use the notation
\begin{equation*}
\widetilde \calT^\lambda(\zeta) = \calT_{\widetilde \varphi^\lambda | \varphi^{\omega_1}}(\zeta), \qquad \calT^\lambda(\zeta) = \calT_{\varphi^\lambda | \varphi^{\omega_1}}(\zeta).
\end{equation*}

\subsection{\texorpdfstring{$L$-operators and $Q$-operators}{L-operators and Q-operators}}

Now let $\varphi$ be a representation of $\uqlbp$ on a vector space $W$, which cannot be extended to a representation of $\uqlg$, and $\psi$ a representation of $\uqlg$ on a vector space $U$. In this case the integrability object $\underset{n}{\bm X}{}_{\varphi | \psi}(\zeta)$ defined by equation (\ref{xfpz}) is called an $L$-operator and denoted $\underset{n}{\bm L}{}_{\varphi | \psi}(\zeta)$.

The companion of an $L$-operator $\underset{n}{\bm L}{}_{\varphi | \psi}(\zeta)$ of type $Y$ is called a $Q$-operator and is denoted $\underset{n}{\bm Q}{}_{\varphi | \psi}(\zeta)$. Explicitly we have
\begin{equation*}
\underset{n}{\bm Q}{}_{\varphi | \psi}(\zeta) = \big( \tr_{\End(W)} \otimes \id_{\End(W^{\otimes n})} \big) \big( \underset{n}{\bm L}{}_{\varphi | \psi}(\zeta) (\varphi_\zeta(t) \otimes 1_{\End(W^{\otimes n})}) \big),
\end{equation*}
where $t$ is a twisting element, which we define by equation (\ref{tqfh}).

Let $\varphi_1$ and $\varphi_2$ be representations of $\uqlbp$, which cannot be extended to representations of $\uqlg$, and $\psi$ a representation of $\uqlg$. It follows from equation (\ref{yyyy}) that
\begin{equation*}
\underset{n}{\bm Q}{}_{\varphi_1 | \psi}(\zeta_1) \, \underset{n}{\bm Q}{}_{\varphi_2 | \psi}(\zeta_2) = \underset{n}{\bm Q}{}_{\varphi_2 | \psi}(\zeta_2) \, \underset{n}{\bm Q}{}_{\varphi_1 | \psi}(\zeta_1),
\end{equation*}
for any $\zeta_1, \zeta_2 \in \bbC^\times$, and similarly for universal transfer operators,
\begin{equation*}
\calQ_{\varphi_1}(\zeta_1) \calQ_{\varphi_2}(\zeta_2) = \calQ_{\varphi_2}(\zeta_2) \calQ_{\varphi_1}(\zeta_1).
\end{equation*}

In the present paper we work with with the $Q$-operators defined as
\begin{equation*}
\underset{n}{\bm Q}{}'_a(\zeta) = \underset{n}{\bm Q}{}_{\theta_a | \varphi^{\omega_1}}(\zeta), \qquad a = 1, \ldots, l + 1,
\end{equation*}
where the representations $\theta_a$ are defined by equation (\ref{thetaa}). The corresponding universal $Q$-operators are denoted as
\begin{equation*}
\calQ'_a(\zeta) = \calQ_{\theta_a}(\zeta).
\end{equation*}
We use a prime to indicate that we redefine these operators below.

\section{Factorization of transfer operators}

\subsection{\texorpdfstring{$l = 1$}{l = 1}}

Consider the tensor product $(W_1)_{\zeta_1} \otimes_\Delta (W_2)_{\zeta_2}$ of two oscillator modules, and introduce an independent labeling for oscillators, using $n_{1 1}$ and $n_{2 1}$ for the first and second factors, respectively. In other words, we label the basis vectors of the tensor product by the 2-tuple of nonnegative integers $(n_{1 1}, \, n_{2 1})$ denoted as $\bm n$. We hope that this slight abuse of notation will not lead to confusion. Compare the product of the highest $\ell$-weights of the modules $(W_1)_{\zeta_1}$ and $ (W_2)_{\zeta_2}$ with the highest $\ell$-weight of the evaluation $\uqlslii$-module $(\widetilde V^\mu)_\zeta$.\footnote{Here and henceforth we treat any $\uqlsllpo$-module as the corresponding $\uqlbp$-module.} The explicit form of these $\ell$-weights is 
\begin{align}
& \bm \Psi_{1, \, \bm 0}(\zeta_1) \bm \Psi_{2, \, \bm 0}(\zeta_2) = \big( - 2 \, \omega_1, \, \big( (1 - q \zeta_2^s u)(1 - q^{- 1} \zeta_1^s u)^{-1} \big) \big), \label{phizphiz} \\[.5em]
& \bm \Lambda^\mu_{\bm 0}(\zeta) = \big( (\mu_1 - \mu_2) \, \omega_1, \, (1 - q^{2 \mu_2} \zeta^s u)(1 - q^{2 \mu_1} \zeta^s u)^{-1} \big),
\end{align}
see sections \ref{s:or} and \ref{s:lwer}. Putting in (\ref{phizphiz})
\begin{equation*}
\zeta_1 = \zeta^\mu_1 = q^{2 (\mu_1 + 1/2) / s} \zeta, \qquad \zeta_2 = \zeta^\mu_2 = q^{2 (\mu_2 - 1/2) / s} \zeta,
\end{equation*}
we obtain
\begin{equation*}
\bm \Psi_{1, \, \bm 0}(\zeta_1) \bm \Psi_{2, \, \bm 0}(\zeta_2) = \big( - 2 \, \omega_1, \, (1 - q^{2 \mu_2} \zeta^s u)(1 - q^{2 \mu_1} \zeta^s u)^{-1} \big).
\end{equation*}
Thus, the product of the highest $\ell$-weights of the modules $(W_1)_{\zeta^\mu_1}$ and $(W_2)_{\zeta^\mu_2}$ coincide with the highest $\ell$-weight of the evaluation $\uqlslii$-module $(\widetilde V^\mu)_\zeta$ shifted by $(- 2 - \mu_1 + \mu_2) \, \omega_1$.

Denote the $\ell$-weights of the module 
\begin{equation*}
W^\mu(\zeta) = (W_1)_{\zeta^\mu_1} \otimes_\Delta (W_2)_{\zeta^\mu_2}
\end{equation*}
as $\bm \Xi^\mu_{\bm n}(\zeta) = (\xi^\mu_{\bm n}, \, \bm \Xi^{\mu +}_{\bm n}(\zeta))$. Using expressions from section \ref{s:or}, we obtain that
\begin{equation*}
\xi^\mu_{\bm n} = (- 2 - 2 n_{1 1} - 2 n_{2 1}) \, \omega_1
\end{equation*}
and $\bm \Xi^{\mu +}_{\bm n}(\zeta) = \big( \Xi^{\mu +}_{\bm n, \, 1}(\zeta, \, u) \big)$, where
\begin{equation*}
\Xi^{\mu +}_{\bm n, \, 1}(\zeta, \, u) = \frac{1 - q^{2 \mu_1 + 2} \zeta^s u}{1 - q^{2 \mu_1 - 2 n_{1 1} + 2)} \zeta^s u} \, \frac{1 - q^{2 \mu_2} \zeta^s u}{1 - q^{2 \mu_1 - 2 n_{1 1}} \zeta^s u}.
\end{equation*}
It is remarkable that $\bm \Xi^{\mu +}_{\bm n}(\zeta)$ does not depend on $n_{2 1}$. In fact, changing $n_{1 1}$ to $m_{1 2}$, we obtain the component $\bm \Lambda^{\mu +}_{\bm m}(\zeta)$ of the $\ell$-weight $\bm \Lambda^\mu_{\bm m}(\zeta)$ of the evaluation $\uqlslii$-module $(\widetilde V^\mu)_\zeta$. Analyzing the expression for the component $\xi^\mu_{\bm n}$ and identifying $n_{1 1}$ with $m_{1 2}$, we see that for any fixed $n_{2 1}$ the $\ell$-weights $\bm \Xi^{\mu +}_{\bm n}(\zeta)$ coincide with the corresponding $\ell$-weights of the module $(\widetilde V^\mu)_\zeta$ shifted by\footnote{We denote $\bm n' = (n_{2 1})$. It may seem that we are using unnecessarily cumbersome notation. This is justified by the fact that we are going to use the same notation for all values of $l$.}
\begin{equation*}
\delta^\mu_{\bm n'} = (- 2 - \mu_1 + \mu_2 - 2 n_{2 1}) \, \omega_1,
\end{equation*}
see equations (\ref{psimu}) and (\ref{psimup}). Thus we have the following relation satisfied by $q$-characters
\begin{equation*}
\chi_q \big (W^\mu(\zeta) \big) = \sum_{\bm n'} \chi_q \big( \widetilde V^\mu_\zeta [\delta_{\bm n'}] \big),
\end{equation*}
which is equivalent to the relation in the Grothendieck ring
\begin{equation*}
\big \langle \, W^\mu(\zeta) \, \big \rangle = \sum_{\bm n'} \big \langle \, \widetilde V^\mu_\zeta [\delta_{\bm n'}] \, \big \rangle = \big \langle \bigoplus_{\bm n'} \widetilde V^\mu_\zeta [\delta_{\bm n'}] \, \big \rangle.
\end{equation*}
Here the summation over $\bm n'$ means the summation over $n_{2 1}$ from $0$ to $\infty$. It follows that
\begin{equation*}
\calY_{W^\mu(\zeta)} = \sum_{\bm n'} \calY_{ \widetilde V^\mu_\zeta [\delta_{\bm n'}]},
\end{equation*}
and, using (\ref{yfsp}), we obtain
\begin{equation*}
\calY_{W^\mu(\zeta)} = q^{(\mu_1 - \mu_2 + 2) \svee[1]{\tilde \alpha} / 2} \sum_{n_{1 2} = 0}^\infty q^{n_{1 2} \, \svee[1]{\tilde \alpha}} \, \widetilde \calT^\mu(\zeta) = q^{(\mu_1 - \mu_2) \svee[1]{\tilde \alpha} / 2} \frac{q^{\svee[1]{\tilde \alpha}}}{1 - q^{\svee[1]{\tilde \alpha}}} \, \widetilde \calT^\mu(\zeta).
\end{equation*}
From the other hand, taking into account (\ref{yffp}), we see that
\begin{equation*}
\calY_{W^\mu(\zeta)} = \calY_{\theta_1}(\zeta^\mu_1) \calY_{\theta_2}(\zeta^\mu_2) = \calQ'_1(\zeta^\mu_1) \calQ'_2(\zeta^\mu_2),
\end{equation*}
and come to the following factorization formula
\begin{equation}
q^{(\mu_1 - \mu_2) \svee[1]{\tilde \alpha} / 2} \frac{q^{\svee[1]{\tilde \alpha}}}{1 - q^{\svee[1]{\tilde \alpha}}} \, \widetilde \calT^\mu(\zeta) = \calQ'_1(\zeta^\mu_1) \calQ'_2(\zeta^\mu_2). \label{tmuqq}
\end{equation}
Introducing the following new universal $Q$-operators
\begin{equation*}
\calQ_1(\zeta) = \zeta^{-\svee[1]{\tilde \alpha} s / 4} \calQ'_1(\zeta), \qquad \calQ_2(\zeta) = \zeta^{\svee[1]{\tilde \alpha} s / 4} \calQ'_2(\zeta),
\end{equation*}
we rewrite the factorization formula as
\begin{equation}
\calC_1 \, \widetilde \calT^\mu(\zeta) = \calQ_1(\zeta^\mu_1) \calQ_2(\zeta^\mu_2), \label{ctqq}
\end{equation}
where
\begin{equation*}
\calC_1 = \frac{q^{\svee[1]{\tilde \alpha} / 2}}{1 - q^{\svee[1]{\tilde \alpha}}}.
\end{equation*}
The advantage of this formula over (\ref{tmuqq}) is that the factor $\calC_1$ does not depend on $\mu$, and this is necessary to obtain the determinant formula.

In terms of usual integrability objects equation (\ref{ctqq} looks as
\begin{equation*}
\rho_{\widetilde \varphi^\mu | \varphi^{\omega_1}}(\zeta)^n \, \underset{n}{\bm C}{}_1 \, \, \underset{n}{\widetilde{\bm T}}{}^\mu (\zeta) = \rho_{\theta_1 | \varphi^{\omega_1}}(\zeta_1^\mu)^n \rho_{\theta_2 | \varphi^{\omega_1}}(\zeta_2^\mu)^n \,  \underset{n}{\bm Q}{}_1(\zeta_1^\mu) \, \underset{n}{\bm Q}{}_2(\zeta_2^\mu),
\end{equation*}
where
\begin{equation*}
\underset{n}{\bm C}{}_1 = (\underbracket[.5pt]{\varphi^{\omega_1} \otimes \cdots \varphi^{\omega_1}}_n) (\calC_1).
\end{equation*}
If we choose
\begin{equation*}
\rho_{\widetilde \varphi^\mu | \varphi^{\omega_1}}(\zeta) = \rho_{\theta_1 | \varphi^{\omega_1}}(\zeta_1^\mu) \rho_{\theta_2 | \varphi^{\omega_1}}(\zeta_2^\mu),
\end{equation*}
we come to simpler expression
\begin{equation*}
\underset{n}{\bm C}{}_1 \, \underset{n}{\widetilde{\bm T}}{}^\mu (\zeta) = \underset{n}{\bm Q}{}_1(\zeta_1^\mu) \, \underset{n}{\bm Q}{}_2(\zeta_2^\mu).
\end{equation*}

\subsection{\texorpdfstring{$l = 2$}{l = 2}}

Consider now the tensor product $(W_1)_{\zeta_1} \otimes_\Delta (W_2)_{\zeta_2} \otimes_\Delta  (W_3)_{\zeta_3}$ of three oscillator modules, and introduce an independent labeling for oscillators, using for it the tuple $\bm n = (n_{1 1}, \, n_{1 2}, \, n_{2 1}, \, n_{2 2}, \, n_{3 1}, \, n_{3 2})$. The product of the highest $\ell$-weights of the modules $(W_1)_{\zeta_1}$, $(W_2)_{\zeta_2}$ and $(W_3)_{\zeta_3}$ is
\begin{multline*}
\bm \Psi_{1, \, \bm 0}(\zeta_1) \bm \Psi_{2, \, \bm 0}(\zeta_2) \bm \Psi_{3, \, \bm 0}(\zeta_3) \\
= \big( - 2 \, \omega_1 - 2 \, \omega_2, \, \big( (1 - \zeta_2^s u)(1 - q^{- 2} \zeta_1^s u)^{-1}, \, (1 - q \zeta_3^s u)(1 - q^{- 1} \zeta_2^s u)^{-1} \big) \big),
\end{multline*}
see section \ref{s:or}, while the highest $\ell$-weight of the evaluation $\uqlsliii$-module $(\widetilde V^\mu)_\zeta$ has the form 
\begin{multline*}
\bm \Lambda^{\mu + }_{\bm 0}(\zeta) = \big( (\mu_1 - \mu_2) \, \omega_1 - (\mu_2 - \mu_3) \, \omega_2, \\ \big( (1 - q^{2 \mu_2} \zeta^s u)(1 - q^{2 \mu_1} \zeta^s u)^{-1}, \, (1 - q^{2 \mu_3 - 1} \zeta^s u)(1 - q^{2 \mu_2 - 1} \zeta^s u)^{-1} \big) \big),
\end{multline*}
see section \ref{s:lwer}. Assuming that 
\begin{equation*}
\zeta_1 = \zeta^\mu_1 = q^{2 (\mu_1 + 1) / s} \zeta, \qquad \zeta_2 = \zeta^\mu_2 = q^{2 \mu_2 / s} \zeta. \qquad \zeta_3 = \zeta^\mu_3 = q^{2 (\mu_3 - 1) / s} \zeta,
\end{equation*}
we see that in this case the product of the highest $\ell$-weights of the modules $(W_1)_{\zeta^\mu_1}$, $(W_2)_{\zeta^\mu_2}$ and $(W_3)_{\zeta^\mu_3}$ coincide with the highest $\ell$-weight of the evaluation $\uqlsliii$-module $(\widetilde V^\mu)_\zeta$ shifted by $(- 2 - \mu_1 + \mu_2) \, \omega_1 + (- 2 - \mu_2 + \mu_3) \, \omega_2$.

Similarly as above, denote the $\ell$-weights of the module 
\begin{equation*}
W^\mu(\zeta) = (W_1)_{\zeta^\mu_1} \otimes_\Delta  (W_2)_{\zeta^\mu_2} \otimes_\Delta  (W_3)_{\zeta^\mu_3},
\end{equation*}
as $\bm \Xi^\mu_{\bm n}(\zeta) = (\xi^\mu_{\bm n}, \, \bm \Xi^{\mu +}_{\bm n}(\zeta))$. Using expressions from section \ref{s:or}, we obtain that
\begin{multline}
\xi^\mu_{\bm n} = (- 2 - 2 n_{1 1} - n_{1 2} + n_{2 1} - 2 n_{2 2} - n_{3 1} + n_{3 2}) \, \omega_1 \\
+  (- 2 + n_{1 1} - n_{1 2} - 2 n_{2 1} + n_{2 2} - n_{3 1} - 2 n_{3 2}) \, \omega_2 \label{ximunb}
\end{multline}
and $\bm \Xi^{\mu +}_{\bm n}(\zeta) = (\Xi^{\mu +}_{\bm n, \, 1}(\zeta, \, u), \, \Xi^{\mu +}_{\bm n, \, 2}(\zeta, \, u))$, where
\begin{align*}
& \Xi^{\mu +}_{\bm n, \, 1}(\zeta, \, u) = \frac{1 - q^{2 \mu_1 - 2 n_{1 2} + 2} \zeta^s u}{1 - q^{2 \mu_1 - 2 n_{1 1} - 2 n_{1 2} + 2} \zeta^s u} \, \frac{1 - q^{2 \mu_2 - 2 n_{2 1}} \zeta^s u}{1 - q^{2 \mu_1 - 2 n_{1 1} - 2 n_{1 2}} \zeta^s u}, \\
& \Xi^{\mu +}_{\bm n, \, 2}(\zeta, \, u) = \frac{1 - q^{2 \mu_1 - 2 n_{1 1} - 2 n_{1 2} + 1} \zeta^s u}{1 - q^{2 \mu_1 - 2 n_{1 2} + 1} \zeta^s u} \, \frac{1 - q^{2 \mu_1 + 3} \zeta^s u}{1 - q^{2 \mu_1 - 2 n_{1 2} + 3} \zeta^s u}   \\*
& \hspace{12em} {} \times \frac{1 - q^{2 \mu_2 + 1} \zeta^s u}{1 - q^{2 \mu_2 - 2 n_{2 1} + 1} \zeta^s u} \, \frac{1 - q^{2 \mu_3 - 1} \zeta^s u}{1 - q^{2 \mu_2 - 2 n_{2 1} - 1} \zeta^s u}.
\end{align*}
We see that now $\bm \Xi^{\mu +}_{\bm n}(\zeta)$ does not depend on $n_{2 2}$, $n_{3 1}$ and $n_{3 2}$. In fact, identifying $n_{1 1}$ with $m_{1 2}$, $n_{1 2}$ with $m_{1 3}$, and $n_{2 1}$ with $m_{2 3}$, we obtain the component $\bm \Lambda^{\mu +}_{\bm m}(\zeta)$ of the $\ell$-weight $\bm \Lambda^\mu_{\bm m}(\zeta)$ of the evaluation $\uqlsliii$-module $\widetilde V^\mu_\zeta$, see equations (\ref{psimup2a}) and (\ref{psimup2b}). Using the above identification, rewrite equation (\ref{ximunb}) as
\begin{multline*}
\xi^\mu_{\bm n} = (\mu_1 - \mu_2 - 2 m_{1 2} - m_{1 3} + m_{2 3}) \, \omega_1 + (\mu_2 - \mu_3 + m_{1 2} - m_{1 3} - 2 m_{2 3}) \omega_2 \\
+ (- 2 - \mu_1 + \mu_2 - 2 n_{2 2} - n_{3 1} + n_{3 2}) \, \omega_1 +  (- 2 - \mu_2 + \mu_3 + n_{2 2} - n_{3 1} - 2 n_{3 2}) \, \omega_2.
\end{multline*}
We see that after the above identification, for any fixed $n_{2 2}$, $n_{3 1}$ and $n_{3 2}$ the $\ell$-weight components $\bm \Xi^{\mu +}_{\bm n}(\zeta)$ coincide with the corresponding $\ell$-weight components of the evaluation $\uqlsliii$-module $(\widetilde V^\mu)_\zeta$ shifted by\footnote{Now we denote $\bm n' = (n_{2 2}, \, n_{3 1}, \, n_{3 2})$.}
\begin{equation*}
\delta_{\bm n'} = (- 2 - \mu_1 + \mu_2 - 2 n_{2 2} - n_{3 1} + n_{3 2}) \, \omega_1 + (- 2 - \mu_2 + \mu_3 + n_{2 2} - n_{3 1} - 2 n_{3 2}) \, \omega_2,
\end{equation*}
see equation (\ref{psimu2}). Thus, we have the relation
\begin{equation*}
\chi_q(W^\mu(\zeta)) = \sum_{\bm n'} \chi_q(\widetilde V^\mu_\zeta [\delta_{\bm n'}]),
\end{equation*}
 satisfied by $q$-characters, which is equivalent to the relation in the Grothendieck ring
\begin{equation*}
\big \langle \, W^\mu(\zeta) \, \big \rangle = \sum_{\bm n'} \big \langle \, \widetilde V^\mu_\zeta [\delta_{\bm n'}] \, \big \rangle = \big \langle \bigoplus_{\bm n'} \widetilde V^\mu_\zeta [\delta_{\bm n'}] \, \big \rangle.
\end{equation*}
Here the summation over $\bm n'$ means the summation over $n_{2 2}$, $n_{3 1}$, and $n_{3 2}$ from $0$ to $\infty$. It follows that
\begin{equation*}
\calY_{W^\mu(\zeta)} = \sum_{\bm n'} \calY_{ \widetilde V^\mu_\zeta [\delta_{\bm n'}]},
\end{equation*}
and, using (\ref{yfsp}), we obtain
\begin{equation*}
\calY_{W^\mu(\zeta)} =  q^{((2 \mu_1 + \mu_2 - \mu_3) \svee[1]{\tilde \alpha} + (\mu_1 + \mu_2 - 2 \mu_3) \svee[2]{\tilde \alpha})/3} \frac{q^{\svee[1]{\tilde \alpha}}}{1 - q^{\svee[1]{\tilde \alpha}}} \, \frac{q^{\svee[1]{\tilde \alpha} + \svee[2]{\tilde \alpha}}}{1 - q^{\svee[1]{\tilde \alpha} + \svee[2]{\tilde \alpha}}} \, \frac{q^{\svee[2]{\tilde \alpha}}}{1 - q^{\svee[2]{\tilde \alpha}}} \, \tilde \calT^\mu(\zeta).
\end{equation*}
From the other hand, taking again into account equation (\ref{yffp}), we see that
\begin{equation*}
\calY_{W^\mu(\zeta)} = \calY_{\theta_1}(\zeta^\mu_1) \calY_{\theta_2}(\zeta^\mu_2) \calY_{\theta_3}(\zeta^\mu_3) = \calQ'_1(\zeta^\mu_1) \calQ'_2(\zeta^\mu_2) \calQ'_3(\zeta^\mu_3),
\end{equation*}
and come to the factorization formula
\begin{equation*}
q^{((2 \mu_1 + \mu_2 - \mu_3) \svee[1]{\tilde \alpha} + (\mu_1 + \mu_2 - 2 \mu_3) \svee[2]{\tilde \alpha})/3} \frac{q^{\svee[1]{\tilde \alpha}}}{1 - q^{\svee[1]{\tilde \alpha}}} \, \frac{q^{\svee[1]{\tilde \alpha} + \svee[2]{\tilde \alpha}}}{1 - q^{\svee[1]{\tilde \alpha} + \svee[2]{\tilde \alpha}}} \, \frac{q^{\svee[2]{\tilde \alpha}}}{1 - q^{\svee[2]{\tilde \alpha}}} \, \widetilde \calT^\mu(\zeta) = \calQ'_1(\zeta^\mu_1) \calQ'_2(\zeta^\mu_2) \calQ'_3(\zeta^\mu_3).
\end{equation*}
Introducing new universal $Q$-operators
\begin{equation*}
\calQ_1(\zeta) = \zeta^{- (2 \svee[1]{\tilde \alpha} + \svee[2]{\tilde \alpha}) s / 6} \calQ'_1(\zeta), \quad \calQ_2(\zeta) = \zeta^{(\svee[1]{\tilde \alpha} - \svee[2]{\tilde \alpha}) s / 6} \calQ'_2(\zeta), \quad \calQ_3(\zeta) = \zeta^{(\svee[1]{\tilde \alpha} + 2 \svee[2]{\tilde \alpha}) s / 6} \calQ'_3(\zeta),
\end{equation*}
we obtain the factorization formula
\begin{equation}
\calC_2 \, \widetilde \calT^\mu(\zeta) = \calQ_1(\zeta^\mu_1) \calQ_2(\zeta^\mu_2) \calQ_3(\zeta^\mu_3), \label{ctqqq}
\end{equation}
where
\begin{equation*}
\calC_2 = \frac{q^{\svee[1]{\tilde \alpha}/2}}{1 - q^{\svee[1]{\tilde \alpha}}} \, \frac{q^{(\svee[1]{\tilde \alpha} + \svee[2]{\tilde \alpha})/2}}{1 - q^{\svee[1]{\tilde \alpha} + \svee[2]{\tilde \alpha}}} \, \frac{q^{\svee[2]{\tilde \alpha}/2}}{1 - q^{\svee[2]{\tilde \alpha}}}.
\end{equation*}
We again come to a factorization formula with the coefficient $\calC_2$ which does not depend on~$\mu$.

Choosing the appropriate normalization, we write (\ref{ctqqq}) as
\begin{equation*}
\underset{n}{\bm C}{}_2 \, \underset{n}{\widetilde{\bm T}}{}^\mu (\zeta) = \underset{n}{\bm Q}{}_1(\zeta_1^\mu) \, \underset{n}{\bm Q}{}_2(\zeta_2^\mu) \, \underset{n}{\bm Q}{}_3(\zeta_3^\mu).
\end{equation*}

\subsection{Determinant formula}

The Weyl group $W$ of the root system of $\gllpo$ is isomorphic to the symmetric group $\rmS_{l + 1}$. It is generated by simple reflections $r_i$, $i = 1, \ldots, l$. The minimal number of generators $r_i$ necessary to represent an element $w \in W$ is said to be the length of $w$ and is denoted by $l(w)$. It is assumed that the identity element has the length equal to $0$.

Using the quantum version of the Bernstein--Gelfand--Gelfand resolution for the quantum group $\uqgllpo$, see, for example, the papers \cite{Ros91, Mal92, HecKol07}, we obtain the equation
\begin{equation*}
\calT^\mu(\zeta) = \sum_{w \in W} (-1)^{l(w)} \widetilde \calT^{\, w \cdot \mu}(\zeta) = \sum_{w \in \mathrm S_{l + 1}} \sgn(w) \, \widetilde \calT^{\, w \cdot \mu}(\zeta),
\end{equation*}
where $w \cdot \lambda$ means the affine action of $w$ defined as
\begin{equation*}
w \cdot \mu = w(\mu + \rho) - \rho
\end{equation*}
with
\begin{equation*}
\rho = \frac{1}{2} \sum_{\substack{i, j = 1 \\ i < j}}^{l + 1} \alpha_{i j}.
\end{equation*}
Now, using (\ref{ctqq}) and (\ref{ctqqq}), for $l = 1$ and $l = 2$ we come to the determinant formula
\begin{equation*}
\calC_l \, \calT^\mu(\zeta) = \det \big( \calQ_a(q^{(2 (\mu_b + \rho_b)/s} \zeta) \big)_{a, \, b = 1}^{l + 1}.
\end{equation*}
Assuming the appropriate normalization, in terms of the usual integrability objects we have
\begin{equation*}
\underset{n}{\bm C}{}_l \, \underset{n}{\bm T}{}^\mu(\zeta) = \det \big( \underset{n}{\bm Q}{}_a(q^{(2 (\mu_b + \rho_b)/s} \zeta) \big)_{a, \, b = 1}^{l + 1}.
\end{equation*}

\section{Conclusion}

Analyzing $\ell$--weights of the evaluation and $q$-oscillator representations we have proved the factorization relations for the transfer operators of the quantum integrable systems associated with the quantum loop algebras $\uqlsllpo$ for $l = 1$ and $l = 2$. The results obtained are completely consistent with the results obtained in the papers \cite{BooGoeKluNirRaz14a, NirRaz16a, BooGoeKluNirRaz14b}. More details and a proof of factorization for the case of arbitrary rank will be given in the forthcoming paper.

\vspace{1em}

{\em Acknowledgments.}
This work was supported in part by the RFBR grant \#~20-51-12005. The author is grateful to H.~ Boos, F.~ G\"ohmann, A.~ Kl\"umper, and Kh.~S.~ Nirov, in collaboration with whom some important results, used in this paper, were previously obtained, for useful discussions.

\newcommand{\noopsort}[1]{}
\providecommand{\bysame}{\leavevmode\hbox to3em{\hrulefill}\thinspace}
\providecommand{\href}[2]{#2}


\begin{thebibliography}{10}

\bibitem{BazLukZam96}
V.~V. Bazhanov, S.~L. Lukyanov, and A.~B. Zamolodchikov, \emph{Integrable
  structure of conformal field theory, quantum {K}d{V} theory and thermodynamic
  {B}ethe ansatz}, \href{http://dx.doi.org/10.1007/BF02101898}{Commun. Math.
  Phys.} \textbf{177} (1996), 381--398,
  \href{http://arxiv.org/abs/hep-th/9412229}{{\tt arXiv:hep-th/9412229}}.

\bibitem{BazLukZam97}
V.~V. Bazhanov, S.~L. Lukyanov, and A.~B. Zamolodchikov, \emph{Integrable
  structure of conformal field theory {II}. {Q}-operator and {DDV} equation},
  \href{http://dx.doi.org/10.1007/s002200050240}{Commun. Math. Phys.}
  \textbf{190} (1997), 247--278,
  \href{http://arxiv.org/abs/hep-th/9604044}{{\tt arXiv:hep-th/9604044}}.

\bibitem{BazLukZam99}
V.~V. Bazhanov, S.~L. Lukyanov, and A.~B. Zamolodchikov, \emph{Integrable
  structure of conformal field theory {III}. {T}he {Y}ang--{B}axter relation},
  \href{http://dx.doi.org/10.1007/s002200050531}{Commun. Math. Phys.}
  \textbf{200} (1999), 297--324,
  \href{http://arxiv.org/abs/hep-th/9805008}{{\tt arXiv:hep-th/9805008}}.

\bibitem{KhoTol92}
S.~M. Khoroshkin and V.~N. Tolstoy, \emph{The uniqueness theorem for the
  universal {$R$}-matrix}, \href{http://dx.doi.org/10.1007/BF00402899}{Lett.
  Math. Phys.} \textbf{24} (1992), 231--244.

\bibitem{LevSoiStu93}
S.~Levendorskii, Ya. Soibelman, and V.~Stukopin, \emph{The quantum {W}eyl group
  and the universal quantum {$R$}-matrix for affine {L}ie algebra
  {$A_1^{(1)}$}}, \href{http://dx.doi.org/10.1007/BF00777372}{Lett. Math.
  Phys.} \textbf{27} (1993), 253--264.

\bibitem{ZhaGou94}
Y.-Z. Zhang and M.~D. Gould, \emph{Quantum affine algebras and universal
  ${R}$-matrix with spectral parameter},
  \href{http://dx.doi.org/10.1007/BF00750144}{Lett. Math. Phys.} \textbf{31}
  (1994), 101--110, \href{http://arxiv.org/abs/hep-th/9307007}{{\tt
  arXiv:hep-th/9307007}}.

\bibitem{BraGouZhaDel94}
A.~J. Bracken, M.~D. Gould, Y.-Z. Zhang, and G.~W. Delius, \emph{Infinite
  families of gauge-equivalent {$R$}-matrices and gradations of quantized
  affine algebras}, \href{http://dx.doi.org/10.1142/S0217979294001585}{Int. J.
  Mod. Phys. B} \textbf{8} (1994), 3679--3691,
  \href{http://arxiv.org/abs/hep-th/9310183}{{\tt arXiv:hep-th/9310183}}.

\bibitem{BraGouZha95}
A.~J. Bracken, M.~D. Gould, and Y.-Z. Zhang, \emph{Quantised affine algebras
  and parameter-dependent {$R$}-matrices},
  \href{http://dx.doi.org/10.1017/S0004972700014040}{Bull. Austral. Math. Soc.}
  \textbf{51} (1995), 177--194.

\bibitem{BooGoeKluNirRaz10}
H.~Boos, F.~G{\"o}hmann, A.~Kl{\"u}mper, Kh.~S. Nirov, and A.~V. Razumov,
  \emph{Exercises with the universal {$R$}-matrix},
  \href{http://dx.doi.org/10.1088/1751-8113/43/41/415208}{J. Phys. A: Math.
  Theor.} \textbf{43} (2010), 415208 (35pp),
  \href{http://arxiv.org/abs/1004.5342}{{\tt arXiv:1004.5342 [math-ph]}}.

\bibitem{BooGoeKluNirRaz11}
H.~Boos, F.~G{\"o}hmann, A.~Kl{\"u}mper, Kh.~S. Nirov, and A.~V. Razumov,
  \emph{On the universal ${R}$-matrix for the {I}zergin--{K}orepin model},
  \href{http://dx.doi.org/10.1088/1751-8113/44/35/355202}{J. Phys. A: Math.
  Theor.} \textbf{44} (2011), 355202 (25pp),
  \href{http://arxiv.org/abs/1104.5696}{{\tt arXiv:1104.5696 [math-ph]}}.

\bibitem{BazTsu08}
V.~V. Bazhanov and Z.~Tsuboi, \emph{Baxter's {Q}-operators for supersymmetric
  spin chains}, \href{http://dx.doi.org/10.1016/j.nuclphysb.2008.06.025}{Nucl.
  Phys. B} \textbf{805} (2008), 451--516,
  \href{http://arxiv.org/abs/0805.4274}{{\tt arXiv:0805.4274 [hep-th]}}.

\bibitem{BooGoeKluNirRaz13}
H.~Boos, F.~G{\"o}hmann, A.~Kl{\"u}mper, Kh.~S. Nirov, and A.~V. Razumov,
  \emph{Universal integrability objects},
  \href{http://dx.doi.org/10.1007/s11232-013-0002-8}{Theor. Math. Phys.}
  \textbf{174} (2013), 21--39, \href{http://arxiv.org/abs/1205.4399}{{\tt
  arXiv:1205.4399 [math-ph]}}.

\bibitem{Raz13}
A.~V. Razumov, \emph{Monodromy operators for higher rank},
  \href{http://dx.doi.org/10.1088/1751-8113/46/38/385201}{J. Phys. A: Math.
  Theor.} \textbf{46} (2013), 385201 (24pp),
  \href{http://arxiv.org/abs/1211.3590}{{\tt arXiv:1211.3590 [math.QA]}}.

\bibitem{BooGoeKluNirRaz14a}
H.~Boos, F.~G{\"o}hmann, A.~Kl\"umper, Kh.~S. Nirov, and A.~V. Razumov,
  \emph{Universal ${R}$-matrix and functional relations},
  \href{http://dx.doi.org/10.1142/S0129055X14300052}{Rev. Math. Phys.}
  \textbf{26} (2014), 1430005 (66pp),
  \href{http://arxiv.org/abs/1205.1631}{{\tt arXiv:1205.1631 [math-ph]}}.

\bibitem{BazHibKho02}
V.~V. Bazhanov, A.~N. Hibberd, and S.~M. Khoroshkin, \emph{Integrable structure
  of {$\mathcal W_3$} conformal field theory, quantum {B}oussinesq theory and
  boundary affine {T}oda theory},
  \href{http://dx.doi.org/10.1016/S0550-3213(01)00595-8}{Nucl. Phys. B}
  \textbf{622} (2002), 475--574,
  \href{http://arxiv.org/abs/hep-th/0105177}{{\tt arXiv:hep-th/0105177}}.

\bibitem{Koj08}
T.~Kojima, \emph{Baxter's ${Q}$-operator for the ${W}$-algebra ${W_N}$},
  \href{http://dx.doi.org/10.1088/1751-8113/41/35/355206}{J. Phys. A: Math.
  Theor} \textbf{41} (2008), 355206 (16pp),
  \href{http://arxiv.org/abs/0803.3505}{{\tt arXiv:0803.3505 [nlin.SI]}}.

\bibitem{BooGoeKluNirRaz14b}
H.~Boos, F.~G{\"o}hmann, A.~Kl\"umper, Kh.~S. Nirov, and A.~V. Razumov,
  \emph{Quantum groups and functional relations for higher rank},
  \href{http://dx.doi.org/10.1088/1751-8113/47/27/275201}{J. Phys. A: Math.
  Theor.} \textbf{47} (2014), 275201 (47pp),
  \href{http://arxiv.org/abs/1312.2484}{{\tt arXiv:1312.2484 [math-ph]}}.

\bibitem{NirRaz14}
Kh.~S. Nirov and A.~V. Razumov, \emph{Quantum groups and functional relations
  for lower rank}, \href{http://arxiv.org/abs/1412.7342}{{\tt arXiv:1412.7342
  [math-ph]}}.

\bibitem{KluNirRaz20}
A.~Kl\"umper, Kh.~S. Nirov, and A.~V. Razumov, \emph{Reduced {qKZ} equation:
  general case}, \href{http://dx.doi.org/10.1088/1751-8121/ab3b9e}{J. Phys. A:
  Math. Gen.} \textbf{53} (2020), 015202 (35pp),
  \href{http://arxiv.org/abs/1905.06014}{{\tt arXiv:1905.06014 [math-ph]}}.

\bibitem{Raz20}
A.~V. Razumov, \emph{Reduced {qKZ} equation and genuine {qKZ} equation},
  \href{http://dx.doi.org/10.1088/1751-8121/aba91d}{J. Phys. A: Math. Theor.}
  (2020), 405204 (32pp), \href{http://arxiv.org/abs/2004.02624}{{\tt
  arXiv:2004.02624 [math-ph]}}.

\bibitem{Res83a}
N.~Yu. Reshetikhin, \emph{The functional equation method in the theory of
  exactly soluble quantum systems}, Sov. Phys. JETP \textbf{57} (1983),
  691--696.

\bibitem{Res83}
N.~Yu. Reshetikhin, \emph{A method of functional equations in the theory of
  exactly solvable quantum systems},
  \href{http://dx.doi.org/10.1007/BF00400435}{Lett. Math. Phys.} \textbf{7}
  (1983), 205--213.

\bibitem{KulRes86}
P.~P. Kulish and N.~Yu. Reshetikhin, \emph{{GL$_3$}-invariant solutions of the
  {Y}ang--{B}axter equation and associated quantum systems}, J. Sov. Math.
  \textbf{34} (1986), 1948--1971.

\bibitem{BazRes90}
V.~V. Bazhanov and N.~Reshetikhin, \emph{Restricted solid-on-solid models
  connected with simply laced algebras and conformal field theory},
  \href{http://dx.doi.org/10.1088/0305-4470/23/9/012}{J. Phys. A: Math. Gen.}
  \textbf{23} (1990), 1477--1492.

\bibitem{KluPea92}
A.~Kl\"umper and P.~A. Pearce, \emph{Conformal weights of {RSOS} lattice models
  and their fusion hierarchies},
  \href{http://dx.doi.org/10.1016/0378-4371(92)90149-K}{Physica A} \textbf{183}
  (1992), 304--350.

\bibitem{KunNakSuz94}
A.~Kuniba, T.~Nakanishi, and J.~Suzuki, \emph{Functional relations in solvable
  lattice models. {I}. {F}unctional relations and representation theory},
  \href{http://dx.doi.org/10.1142/S0217751X94002119}{Int. J. Mod. Phys. A}
  \textbf{9} (1994), 5215--5266,
  \href{http://arxiv.org/abs/hep-th/9309137}{{\tt arXiv:hep-th/9309137}}.

\bibitem{KunNakSuz11}
A.~Kuniba, T.~Nakanishi, and J.~Suzuki, \emph{${T}$-systems and ${Y}$-systems
  in integrable systems},
  \href{http://dx.doi.org/10.1088/1751-8113/44/10/103001}{J. Phys A: Math.
  Theor.} (2011), 103001 (146pp), \href{http://arxiv.org/abs/1010.1344}{{\tt
  arXiv:1010.1344 [hep-th]}}.

\bibitem{NirRaz16a}
Kh.~S. Nirov and A.~V. Razumov, \emph{Quantum groups and functional relations
  for lower rank}, \href{http://dx.doi.org/10.1016/j.geomphys.2016.10.014}{J.
  Geom. Phys.} \textbf{112} (2017), 1--28,
  \href{http://arxiv.org/abs/1412.7342}{{\tt arXiv:1412.7342 [math-ph]}}.

\bibitem{FreHer15}
E.~Frenkel and D.~Hernandez, \emph{Baxter's relations and spectra of quantum
  integrable models}, Duke Math. J. \textbf{164} (2015), 2407--2460,
  \href{http://arxiv.org/abs/1308.3444}{{\tt arXiv:1308.3444 [math.QA]}}.

\bibitem{HerJim12}
D.~Hernandez and M.~Jimbo, \emph{Asymptotic representations and {D}rinfeld
  rational fractions}, \href{http://dx.doi.org/10.1112/S0010437X12000267}{Comp.
  Math.} \textbf{148} (2012), 1593--1623,
  \href{http://arxiv.org/abs/1104.1891}{{\tt arXiv:1104.1891 [math.QA]}}.

\bibitem{Jim86a}
M.~Jimbo, \emph{A $q$-analogue of {$\mathrm U(\mathfrak{gl}(N + 1))$}, {H}ecke
  algebra, and the {Y}ang--{B}axter equation},
  \href{http://dx.doi.org/10.1007/BF00400222}{Lett. Math. Phys.} \textbf{11}
  (1986), 247--252.

\bibitem{LezSav74}
A.~N. Leznov and M.~V. Saveliev, \emph{A parametrization of compact groups},
  \href{http://dx.doi.org/10.1007/BF01075497}{Funct. Anal. Appl.} \textbf{8}
  (1974), 347--348.

\bibitem{AshSmiTol79}
R.~M. Asherova, Yu.~F. Smirnov, and V.~N. Tolstoy, \emph{Description of a class
  of projection operators for semisimple complex {L}ie algebras},
  \href{http://dx.doi.org/10.1007/BF01140268}{Math. Notes} \textbf{26} (1979),
  499--504.

\bibitem{Tol89}
V.~N. Tolstoy, \emph{Extremal projections for contragredient {L}ie algebras and
  superalgebras of finite growth},
  \href{http://dx.doi.org/10.1070/RM1989v044n01ABEH002023}{Russian Math.
  Surveys} \textbf{44} (1989), 257--258.

\bibitem{NirRaz17b}
Kh.~S. Nirov and A.~V. Razumov, \emph{Quantum groups, {V}erma modules and
  $q$-oscillators: general linear case},
  \href{http://dx.doi.org/10.1088/1751-8121/aa7808}{J. Phys. A: Math. Theor.}
  \textbf{50} (2017), 305201 (19pp),
  \href{http://arxiv.org/abs/1610.02901}{{\tt arXiv:1610.02901 [math-ph]}}.

\bibitem{Yam89}
H.~Yamane, \emph{A {P}oincar\'e--{B}irkhoff--{W}itt theorem for quantized
  universal enveloping algebras of type {$A_N$}},
  \href{http://dx.doi.org/10.2977/prims/1195173355}{Publ. RIMS. Kyoto Univ.}
  \textbf{25} (1989), 503--520.

\bibitem{Ser01}
J.-P. Serre, \emph{Complex semisimple {Lie} algebras}, Springer Monographs in
  Mathematics, Springer, Berlin, 2001.

\bibitem{Hum80}
J.~E. Humphreys, \emph{Introduction to {L}ie algebras and representation
  theory}, Springer, New York, 1980.

\bibitem{Kac90}
V.~Kac, \emph{Infinite dimensional {L}ie algebras}, Cambridge University Press,
  Cambridge, 1990.

\bibitem{TolKho92}
V.~N. Tolstoy and S.~M. Khoroshkin, \emph{The universal {$R$}-matrix for
  quantum utwisted affine {L}ie algebras},
  \href{http://dx.doi.org/10.1007/BF01077085}{Funct. Anal. Appl.} \textbf{26}
  (1992), 69--71.

\bibitem{KhoTol93}
S.~M. Khoroshkin and V.~N. Tolstoy, \emph{On {D}rinfeld's realization of
  quantum affine algebras},
  \href{http://dx.doi.org/10.1016/0393-0440(93)90070-U}{J. Geom. Phys.}
  \textbf{11} (1993), 445--452.

\bibitem{KhoTol94}
S.~Khoroshkin and V.~N. Tolstoy, \emph{Twisting of quantum (super)algebras.
  {C}onnection of {D}rinfeld's and {C}artan-{W}eyl realizations for quantum
  affine algebras}, \href{http://arxiv.org/abs/hep-th/9404036}{{\tt
  arXiv:hep-th/9404036}}.

\bibitem{Bec94a}
J.~Beck, \emph{Convex bases of {PBW} type for quantum affine algebras},
  \href{http://dx.doi.org/10.1007/BF02099742}{Commun. Math. Phys.} \textbf{165}
  (1994), 193--199, \href{http://arxiv.org/abs/hep-th/9407003}{{\tt
  arXiv:hep-th/9407003}}.

\bibitem{Dam98}
I.~Damiani, \emph{La {$R$}-matrice pour les alg{\`e}bres quantiques de type
  affine non tordu},
  \href{http://dx.doi.org/10.1016/S0012-9593(98)80104-3}{Ann. Sci. {\'E}cole
  Norm. Sup.} \textbf{31} (1998), 493--523.

\bibitem{Dri87}
V.~G. Drinfeld, \emph{Quantum groups}, Proceedings of the International
  Congress of Mathematicians, Berkeley, 1986 (A.~E. Gleason, ed.), vol.~1,
  American Mathematical Society, Providence, 1987, pp.~798--820.

\bibitem{Dri88}
V.~G. Drinfeld, \emph{A new realization of {Y}angians and quantized affine
  algebras}, Soviet Math. Dokl. \textbf{36} (1988), 212--216.

\bibitem{Lax07}
P.~D. Lax, \emph{Linear algebra and its applications}, Wiley, 2007.

\bibitem{BooGoeKluNirRaz16}
H.~Boos, F.~G\"ohmann, A.~Kl\"umper, Kh.~S. Nirov, and A.~V. Razumov,
  \emph{Oscillator versus prefundamental representations},
  \href{http://dx.doi.org/10.1063/1.4966925}{J. Math. Phys.} \textbf{57}
  (2016), 111702 (23pp), \href{http://arxiv.org/abs/1512.04446}{{\tt
  arXiv:1512.04446 [math-ph]}}.

\bibitem{BooGoeKluNirRaz17b}
H.~Boos, F.~G\"ohmann, A.~Kl\"umper, Kh.~S. Nirov, and A.~V. Razumov,
  \emph{Oscillator versus prefundamental representations {II}. {A}rbitrary
  higher ranks}, \href{http://dx.doi.org/10.1063/1.5001336}{J. Math. Phys.}
  \textbf{58} (2017), 093504 (23pp),
  \href{http://arxiv.org/abs/1701.02627}{{\tt arXiv:1701.02627 [math-ph]}}.

\bibitem{FreRes99}
E.~Frenkel and N.~Reshetikhin, \emph{The $q$-characters of representations of
  quantum affine algebras and deformations of ${W}$-algebras}, Contemp. Math.
  \textbf{248} (1999), 163--205, \href{http://arxiv.org/abs/math/9810055}{{\tt
  arXiv:math/9810055}}.

\bibitem{EtiFreKir98}
P.~Etingof, B.~Frenkel, and A.~A. Kirillov, \emph{Lectures on representation
  theory and {K}nizhnik--{Z}amolodchikov equations}, Mathematical Surveys and
  Monographs, vol.~58, American Mathematical Society, Providence, 1998.

\bibitem{NirRaz19}
Kh.~S. Nirov and A.~V. Razumov, \emph{Vertex models and spin chains in formulas
  and pictures}, \href{http://dx.doi.org/10.3842/SIGMA.2019.068}{SIGMA}
  \textbf{15} (2019), 068 (67pp), \href{http://arxiv.org/abs/1811.09401}{{\tt
  arXiv:1811.09401 [math-ph]}}.

\bibitem{Tan92}
T.~Tanisaki, \emph{{K}illing forms, {H}arish-{C}handra homomorphisms and
  universal ${R}$-matrices for quantum algebras}, Infinite Analysis
  (A.~Tsuchiya, T.~Eguchi, and M.~Jimbo, eds.), Advanced Series in Mathematical
  Physics, vol.~16, World Scientific, Singapore, 1992, pp.~941--962.

\bibitem{ChaPre94}
V.~Chari and A.~Pressley, \emph{A guide to quantum groups}, Cambridge
  University Press, Cambridge, 1994.

\bibitem{Usm94}
R.~A. Usmani, \emph{Inversion of a tridiagonal {J}acobi matrix},
  \href{http://dx.doi.org/10.1016/0024-3795(94)90414-6}{Linear Algebra Appl.}
  \textbf{212/213} (1994), 413--414.

\bibitem{Ros91}
M.~Rosso, \emph{An analogue of {B.G.G}. resolution for the quantum
  ${SL(N)}$-group}, Symplectic geometry and mathematical physics
  (Aix-en-Provence, 1990) (P.~Donato, C.~Duval, J.~Elhadad, and G.~M. Tynman,
  eds.), Progress in Mathematics, vol.~99, Birkh\"auser, Boston, 1991,
  pp.~422--432.

\bibitem{Mal92}
F.~Malikov, \emph{Quantum groups: singular vectors and {BGG} resolution},
  \href{http://dx.doi.org/10.1142/S0217751X92003963}{Int. J. Mod. Phys. A}
  \textbf{7S1B} (1992), 623--643.

\bibitem{HecKol07}
I.~Heckenberger and S.~Kolb, \emph{On the {B}erstein--{G}elfand--{G}elfand
  resolution for {K}ac--{M}oody algebras and quantized enveloping algebras},
  \href{http://dx.doi.org/10.1007/s00031-007-0059-2}{Transform. Groups}
  \textbf{12} (2007), 647--655, \href{http://arxiv.org/abs/math/0605460}{{\tt
  arXiv:math/0605460}}.

\end{thebibliography}
\end{document}